# CENTRAL LIMIT THEOREMS FOR POISSON HYPERPLANE TESSELLATIONS[1]

By Lothar Heinrich, Hendrik Schmidt and Volker Schmidt

*University of Augsburg, University of Ulm and University of Ulm*

We derive a central limit theorem for the number of vertices of convex polytopes induced by stationary Poisson hyperplane processes in $\mathbb{R}^d$. This result generalizes an earlier one proved by Paroux [*Adv. in Appl. Probab.* **30** (1998) 640–656] for intersection points of motion-invariant Poisson line processes in $\mathbb{R}^2$. Our proof is based on Hoeffding's decomposition of $U$-statistics which seems to be more efficient and adequate to tackle the higher-dimensional case than the "method of moments" used in [*Adv. in Appl. Probab.* **30** (1998) 640–656] to treat the case $d = 2$. Moreover, we extend our central limit theorem in several directions. First we consider $k$-flat processes induced by Poisson hyperplane processes in $\mathbb{R}^d$ for $0 \le k \le d-1$. Second we derive (asymptotic) confidence intervals for the intensities of these $k$-flat processes and, third, we prove multivariate central limit theorems for the $d$-dimensional joint vectors of numbers of $k$-flats and their $k$-volumes, respectively, in an increasing spherical region.

**1. Introduction.** Central limit theorems (briefly CLTs) for models of stochastic geometry have been considered in various papers. For example, [1] and [24] investigate CLTs for Poisson–Voronoi and Poisson line tessellations in the Euclidean plane, respectively. More general CLTs for Poisson–Voronoi tessellations in the $d$-dimensional Euclidean space $\mathbb{R}^d$ have been established in [12] and [25]. In [9], normal approximations are given for some mean-value estimates of absolutely regular ($\beta$-mixing) tessellations. A CLT for stationary tessellations with random inner cell structures has been derived in [13]. Furthermore, CLTs and related asymptotic properties for the empirical volume fraction of stationary random sets in $\mathbb{R}^d$ are examined in [2, 5, 16].

Received February 2005; revised November 2005.

[1]Supported by France Télécom Grant 42 36 68 97.

A CLT for estimators of surface area densities in the Boolean model has been proved in [18], while in [11] CLTs for a more general class of random measures associated with absolutely regular germ–grain models have been proved. In [19] (and references therein), the reader can find a lot of further CLTs for empirical characteristics of Boolean models. Consistency properties and asymptotic normality of joint estimators for the whole vector of specific intrinsic volumes of stationary random sets in $\mathbb{R}^d$ have been derived in [23] and [28], while uniformly best unbiased estimators for the intensity of stationary flat processes have been considered in [27].

In the present paper we prove CLTs for the number of vertices and the number, as well as the volume, of $k$-flats ($k=1,\ldots,d-1$) induced by intersections of $d$ (resp. $d-k$) hyperplanes of stationary, not necessarily isotropic Poisson hyperplane processes in $\mathbb{R}^d$; see Sections 3 and 4. More precisely, we count the number of vertices lying in the $d$-dimensional ball $B_r^d$ with radius $r > 0$ and centered at the origin $o \in \mathbb{R}^d$. In addition, we both count the induced $k$-flats hitting $B_r^d$ and measure their total $k$-dimensional Lebesgue volume in $B_r^d$ for any $k = 0, \ldots, d-1$ and we study their joint behavior when the radius $r$ tends to infinity.

Noting that the number of vertices contained in $B_r^d$ can be expressed in the form of a multiple random sum running over all $d$-tuples of distinct hyperplanes which have a common point in $B_r^d$ (see Chapter 6 in [17]), we use Hoeffding's decomposition of $U$-statistics with a Poisson distributed number of random variables. Hence, asymptotic normality of the number of vertices is obtained by proving a CLT for a Poisson distributed number of independent and identically distributed (i.i.d.) random variables; see Theorem 3.1. Using a similar representation as multiple random sum for the number of $k$-flats hitting $B_r^d$ and for the total $k$-volume of the $k$-flats in $B_r^d$, we generalize the latter CLT in Section 4; see Theorem 4.1. Based on these CLTs, we obtain asymptotic confidence intervals for the intensities of the induced $k$-flat processes and, quite naturally, are able to consider the case of multidimensional CLTs.

We should mention that the normalization in our CLTs is, up to certain constants, with respect to the $d$-dimensional volume of $B_r^d$ raised to the power $1 - 1/(2d)$. We may interpret this as an expression of long-range dependences generated by the hyperplanes themselves. Furthermore, the choice of spherical sampling regions simplifies the proofs considerably, however, most of the results remain valid for more general families of increasing convex sampling windows. If, additionally, isotropy is assumed and no restriction is imposed on the orientation vectors of the intersecting hyperplanes, this allows to determine centering and normalizing constants in the CLTs (i.e., intensities and asymptotic variances) explicitly. Moreover, the results of the present paper, together with Lemma 4.1 in [13], which states that the influence of cells hitting the boundary of $B_r^d$ is asymptotically negligible as $r \to \infty$, it is possible to derive CLTs for $k$-facets ($k = 1, \ldots, d$) of Poisson hyperplane tessellations.

In Section 5 we reformulate Theorem 3.1 in the particular case $d = 2$ and compare this CLT with a related result obtained by Paroux [24] for planar Poisson line processes. Applying again Hoeffding's CLT for $U$-statistics (with random normalization), we obtain a considerably simple proof of the CLT derived in [24] by the "method of moments."

Applications for our results arise in stochastic–geometric network modeling, both in macroscopic settings like in telecommunication (see, e.g., [8]) and in microscopic settings like in cell biology (see, e.g., [4]). In particular, in Section 4.2 we show how our central limit theorems and especially our (asymptotic) confidence intervals can be applied in the framework of the so-called stochastic subscriber line model (SSLM) for telecommunication networks in urban environments. The SSLM is used in the context of strategic network planning and network analysis as a flexible model depending only on a limited number of parameters; see [7]. Figure 1 shows a realization of the SSLM in the case where a Poisson line process is used to model the underlying road system and where two types of network components are placed onto the lines. Besides tessellations induced by Poisson line processes, the class of Voronoi type tesselations is also used in the SSLM, for example, in order to model serving zones; see Figure 1. Therefore, we briefly discuss CLTs for Poisson–Voronoi tessellations in Section 6 which recently have been obtained in [12]; see also [25].

**2. Preliminaries.** In this section the basic notation used in the present paper is introduced and a brief account of some relevant notions of stochastic geometry is given. For a detailed discussion of the subject, the reader is referred to the literature, for example, [30] and [32]. Further background about random tessellations, flat and hyperplane processes can be found, for example, in [20] and [22].

Throughout, let $[\Omega, \sigma(\Omega), \mathbb{P}]$ be a common probability space on which all random objects are defined in the present paper. Let $\langle x, y \rangle = \sum_{k=1}^d x_k y_k$ denote the scalar product of the coordinate vectors $x = (x_1, \ldots, x_d)^\top$ and $y = (y_1, \ldots, y_d)^\top$ in $\mathbb{R}^d$. By means of the Euclidean norm $\|\cdot\| = \sqrt{\langle \cdot, \cdot \rangle}$, we may define the ball $B_r^d = \{x \in \mathbb{R}^d : \|x\| \leq r\}$ centered at the origin and the unit sphere $S^{d-1} = \{x \in \mathbb{R}^d : \|x\| = 1\}$ in $\mathbb{R}^d$, respectively. Furthermore, let $S_+^{d-1} = \{(x_1, \ldots, x_d)^\top \in S^{d-1} : x_d \geq 0\}$ be the upper unit hemisphere and let $\nu_k(\cdot)$ denote the Lebesgue measure in $\mathbb{R}^k$; $k = 0, \ldots, d$. This measure will also be used instead of the $k$-dimensional Hausdorff measure in $\mathbb{R}^d$ for $k = 0, \ldots, d-1$. As usual, $\nu_0(\cdot)$ coincides with the counting measure, that

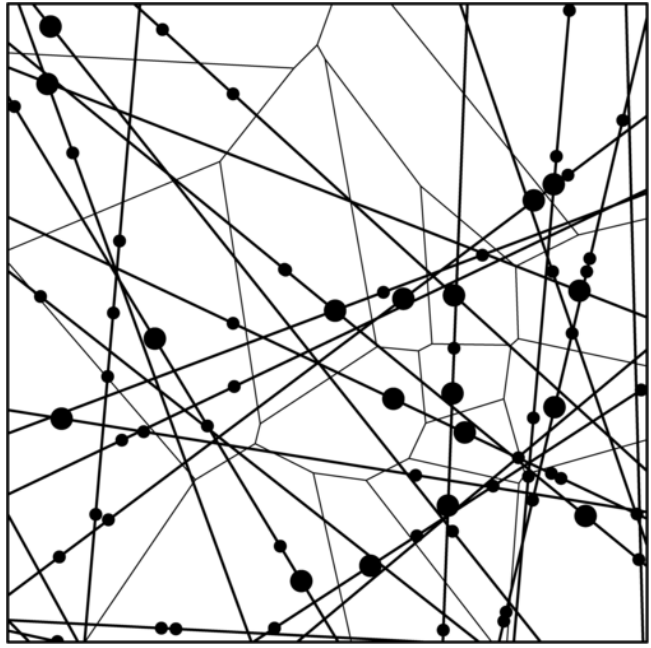

Fig. 1. *Realization of the stochastic subscriber line model.*

is, $\nu_0(B) = \#B$. For brevity, put

$$\kappa_d = \nu_d(B_1^d) = \frac{\pi^{d/2}}{\Gamma(d/2+1)},$$

where $\Gamma(s) = \int_0^\infty e^{-y} y^{s-1}\, dy$ for $s > 0$.

2.1. *Stationary flat processes.* For each $k \in \{0, \ldots, d-1\}$, let $\mathcal{A}_k^d$ denote the space of all affine $k$-dimensional subspaces in $\mathbb{R}^d$ and $\mathcal{L}_k^d = \{L \in \mathcal{A}_k^d : o \in L\}$. A (translation-invariant, resp. rotation-invariant) point process $\Phi_k : \Omega \to N(\mathcal{F}')$ is said to be a (stationary, resp. isotropic) $k$-*flat process* in $\mathbb{R}^d$ if the intensity measure $\Lambda_k(\cdot) = \mathbb{E}\Phi_k(\cdot)$ of $\Phi_k$ satisfies

$$\Lambda_k(\mathcal{F}' \setminus \mathcal{A}_k^d) = 0,$$

where $\mathcal{F}'$ denotes the space of nonempty closed sets in $\mathbb{R}^d$ and $N(\mathcal{F}')$ is the family of all locally finite counting measures on the Borel $\sigma$-algebra $\mathcal{B}(\mathcal{F}')$. A $(d-1)$-flat process is called a *hyperplane process*.

Provided that the intensity measure $\Lambda_k(\cdot)$ of a stationary $k$-flat process $\Phi_k$ is locally finite and different from the zero measure, there exists a finite number $\lambda_k > 0$ (called the *intensity* of $\Phi_k$) and a probability measure $\widetilde{\Theta}_k$ on $\mathcal{B}(\mathcal{L}_k^d)$ (the so-called *orientation distribution* of $\Phi_k$) such that the following

disintegration formula

$$\Lambda_k(B) = \lambda_k \int_{\mathcal{L}_k^d} \int_{L^\perp} \mathbb{1}_B(L+x) \nu_{d-k}(dx) \widetilde{\Theta}_k(dL) \quad (2.1)$$

holds for any $B \in \mathcal{B}(\mathcal{A}_k^d)$. Here, $L^\perp \in \mathcal{L}_{d-k}^d$ denotes the orthogonal complement of $L \in \mathcal{L}_k^d$ and $\mathbb{1}_B(\cdot)$ stands for the indicator function of the set $B$. Formula (2.1) yields a simple interpretation of the intensity $\lambda_k$ as ratio

$$\lambda_k = \frac{\mathbb{E}\Phi_k(\{L \in \mathcal{A}_k^d : L \cap B_r^d \neq \varnothing\})}{\kappa_{d-k} r^{d-k}} \quad \text{for all } r > 0. \quad (2.2)$$

In other words, $\lambda_k \kappa_{d-k}$ is the expected number of $k$-flats hitting the unit ball in $\mathbb{R}^d$. On the other hand, if we use (2.1) and apply Campbell's theorem to the stationary random measure

$$\zeta_{\Phi_k}(\cdot) = \sum_{L \in \mathrm{supp}(\Phi_k)} \nu_k((\cdot) \cap L)$$

on $\mathcal{B}(\mathbb{R}^d)$, where $\mathrm{supp}(\Phi_k) = \{L \in \mathcal{A}_k^d : \Phi_k(L) \geq 1\}$, we get

$$\mathbb{E}\zeta_{\Phi_k}(B) = \lambda_k \nu_d(B) \quad \text{for bounded } B \in \mathcal{B}(\mathbb{R}^d). \quad (2.3)$$

Hence, $\lambda_k$ can be regarded as mean total $k$-volume of all $k$-flats in the unit cube $[0,1)^d$.

In the particular case of a stationary hyperplane process $\Phi$ with intensity $\lambda$, formula (2.1) simplifies since each hyperplane $H(p,v) = \{x : \langle x, v \rangle = p\}$ can be parameterized by its signed perpendicular distance $p \in \mathbb{R}$ from the origin and its orientation vector $v \in S_+^{d-1}$. Using the (spherical) orientation distribution $\Theta : \mathcal{B}(S_+^{d-1}) \to [0,1]$ defined by

$$\Theta_{d-1}(B) = \widetilde{\Theta}(\{H(0,v) \in \mathcal{L}_{d-1}^d : v \in B\}) \quad \text{for } B \in \mathcal{B}(S_+^{d-1}), \quad (2.4)$$

formula (2.1) can be rewritten as

$$\Lambda_{d-1}(B) = \lambda \int_{S_+^{d-1}} \int_{\mathbb{R}} \mathbb{1}_B(H(p,v)) \, dp \, \Theta(dv) \quad \text{for } B \in \mathcal{B}(\mathcal{A}_{d-1}^d).$$

Alternatively, a (spherical) orientation distribution can be introduced as an even (symmetric) probability measure $\Theta^*$ on $\mathcal{B}(S^{d-1})$ which is connected with $\Theta$ by $\Theta^*(B) = \frac{1}{2}(\Theta(B \cap S_+^{d-1}) + \Theta(-B \cap S_+^{d-1}))$ for $B \in \mathcal{B}(S^{d-1})$. The symmetry condition $\Theta^*(B) = \Theta^*(-B)$ expresses the identification of hyperplanes with antipodal orientation vectors. A stationary hyperplane process $\Phi$ is said to be isotropic (or motion-invariant) if $\Theta$ is the uniform distribution on $S_+^{d-1}$ or, equivalently, $\Theta^*$ is the uniform distribution on $S^{d-1}$.

2.2. *Poisson hyperplane processes.* The above parameterization of hyperplanes admits the representation of the (stationary) Poisson hyperplane process $\Phi$ (with intensity $\lambda$ and spherical orientation distribution $\Theta$) as a stationary and independently marked Poisson point process $\Psi = \sum_{i \geq 1} \delta_{(P_i, V_i)}$ on $\mathbb{R}^1$ with intensity $\lambda$ and mark distribution $\Theta$ (on the mark space $S_+^{d-1}$).

Now, let $k \in \{0, \ldots, d-1\}$ be fixed. The intersections of $d-k$ (distinct) hyperplanes $H_1, \ldots, H_{d-k}$ belonging to the support of $\Phi$ induce a stationary $k$-flat process $\Phi_k$ which can be represented as multiple sum, that is,

$$(2.5) \qquad \Phi_k(B) = \frac{1}{(d-k)!} \sum_{H_1, \ldots, H_{d-k} \in \mathrm{supp}(\Phi)} \mathbb{1}_B(H_1 \cap \cdots \cap H_{d-k})$$

for $B \in \mathcal{B}(\mathcal{A}_k^d)$. If the indicator function in the latter expression is replaced by $\nu_k(H_1 \cap \cdots \cap H_{d-k} \cap (\cdot))$, we get a stationary random measure $\zeta_k(\cdot)$ on $\mathcal{B}(\mathbb{R}^d)$. Intensity $\lambda_k$ and orientation distribution $\widetilde{\Theta}_k$ of $\Phi_k$ are given by

$$(2.6) \qquad \begin{aligned} &\lambda_k \widetilde{\Theta}_k(B) \\ &= \frac{\lambda^{d-k}}{(d-k)!} \int_{S_+^{d-1}} \cdots \int_{S_+^{d-1}} \mathbb{1}_B(H(0, v_1) \cap \cdots \cap H(0, v_{d-k})) \\ &\qquad \times \nu_{d-k}(v_1, \ldots, v_{d-k}) \Theta(dv_1) \cdots \Theta(dv_{d-k}), \end{aligned}$$

for $B \in \mathcal{B}(\mathcal{L}_k^d)$; see, for example, Chapter 4 in [30], where $\nu_{d-k}(v_1, \ldots, v_{d-k})$ denotes the $(d-k)$-dimensional volume of the parallelotope spanned by the vectors $v_1, \ldots, v_{d-k} \in S_+^{d-1}$.

The intensity $\lambda_k$ of $\Phi_k$ can be expressed in terms of the $(d-k)$th intrinsic volume of the Steiner convex set (or zonoid) associated with $\Phi$; see, for example, page 161 in [17]. In the isotropic case (i.e., $\Theta$ is the uniform distribution), these formulae reduce to

$$(2.7) \qquad \lambda_k = \binom{d}{k} \frac{\kappa_d}{\kappa_k} \left( \frac{\kappa_{d-1}}{d \kappa_d} \right)^{d-k} \lambda^{d-k} \qquad \text{for } k = 0, 1, \ldots, d-1.$$

2.3. *Hoeffding's decomposition of U-statistics.* The proofs of the central limit theorems we are going to present in Sections 3 and 4 are based on Hoeffding's decomposition of $U$-statistics which we briefly sketch subsequently. A more detailed discussion can be found in [15] and [31].

Let $X_1, X_2, \ldots$ be a sequence of i.i.d. random vectors in $\mathbb{R}^d$ and, for any fixed $m \geq 2$, let $f : \mathbb{R}^{md} \to \mathbb{R}$ be a Borel-measurable symmetric function such that $\mathbb{E}|f(X_1, \ldots, X_m)| < \infty$. A $U$-statistic $U_n^{(m)}(f)$ of order $m$ with *kernel function* $f$ is then defined by

$$(2.8) \qquad U_n^{(m)}(f) = \frac{1}{\binom{n}{m}} \sum_{1 \leq i_1 < \cdots < i_m \leq n} f(X_{i_1}, \ldots, X_{i_m}) \qquad \text{for } n \geq m.$$

Note that $U_n^{(m)}(f)$ is an unbiased estimator for $\mu = \mathbb{E}f(X_1,\ldots,X_m)$. By elementary rearrangements, we may write $U_n^{(m)}(f)$ in the following form (*Hoeffding's decomposition*):

$$U_n^{(m)}(f) - \mu = \frac{m}{n}\sum_{i=1}^{n}(g(X_i) - \mu) + R_n^{(m)}(f), \tag{2.9}$$

where

$$g(x) = \mathbb{E}(f(X_1, X_2, \ldots, X_m)|X_1 = x) = \mathbb{E}f(x, X_2, \ldots, X_m)$$

and

$$R_n^{(m)}(f) = \sum_{k=2}^{m}\binom{m}{k}\binom{n}{k}^{-1}\sum_{1\leq i_1<\cdots<i_k\leq n}G_k(X_{i_1},\ldots,X_{i_k}),$$

with functions $G_k: \mathbb{R}^{kd} \to \mathbb{R}$, $2 \leq k \leq m$, defined by

$$G_k(x_1,\ldots,x_k) = \sum_{j=1}^{k}(-1)^{k-j}\sum_{1\leq i_1<\cdots<i_j\leq k}(\mathbb{E}f(x_{i_1},\ldots,x_{i_j},X_{j+1},\ldots,X_m) - \mu).$$

The crucial outcome of Hoeffding's decomposition (2.9) can be summarized in the estimate

$$\mathbb{E}(R_n^{(m)}(f))^2 \leq \frac{c_m}{n^2}\mathbb{E}f^2(X_1,\ldots,X_m) \qquad \text{for } n \geq m \tag{2.10}$$

and for some constant $c_m < \infty$ only depending on $m$. The latter result, provided that $\mathbb{E}f^2(X_1,\ldots,X_m) < \infty$, immediately leads to Hoeffding's CLT for $U$-statistics (see Chapter 5.2 in [31]), that is,

$$\sqrt{n}(U_n^{(m)}(f) - \mu) \xrightarrow[n\to\infty]{\text{d}} \mathcal{N}(0, m^2(\mathbb{E}g^2(X_1) - \mu^2)), \tag{2.11}$$

where $\xrightarrow{\text{d}}$ denotes convergence in distribution and where $\mathcal{N}(0, m^2(\mathbb{E}g^2(X_1) - \mu^2))$ is a Gaussian mean 0 random variable with variance $m^2(\mathbb{E}g^2(X_1) - \mu^2)$.

**3. Point process of intersection points.** Let $\Psi = \sum_{i\geq 1}\delta_{(P_i,V_i)}$ be the marked-point-process representation of a stationary (not necessarily isotropic) Poisson hyperplane process $\Phi$ with intensity $\lambda > 0$ and a nondegenerate spherical orientation distribution $\Theta$, that is, $\Theta(L \cap S_+^{d-1}) < 1$ for any $L \in \mathcal{L}_{d-1}^d$. This assumption on $\Theta$ ensures that each of the stationary $k$-flat processes $\Phi_k$ generated by $\Phi$ has positive intensity $\lambda_k$ for $k = 0, \ldots, d-1$ and the Poisson hyperplane tessellation induced by $\Phi$ consists of bounded cells; see Chapter 6 in [30].

In this and the next section we derive CLTs for the number $\Phi_k(\{L \in \mathcal{A}_k^d : L \cap B_r^d \neq \varnothing\})$ of $k$-flats ($k = 0, \ldots, d-1$) hitting the ball $B_r^d$, as well as for their total $k$-volume contained in $B_r^d$ when the radius $r$ tends to infinity. In the particular case $k = 0$, the atoms of the point process $\Phi_0$ will be labeled by the $d$ orientation vectors of the intersecting hyperplanes generating the intersection points. More precisely, for any $r > 0$ and $B \in \mathcal{B}((S_+^{d-1})^d)$, we consider the number $\Psi_0(B_r^d \times B)$ of those intersection points $H(P_{i_1}, V_{i_1}) \cap \cdots \cap H(P_{i_d}, V_{i_d})$ in $B_r^d$ for which the corresponding orientation vectors $V_{i_1}, \ldots, V_{i_d}$ satisfy the condition $(V_{(i_1)}, \ldots, V_{(i_d)}) \in B$. Here, $(V_{(i_1)}, \ldots, V_{(i_d)})$ is the re-ordered vector $(V_{i_1}, \ldots, V_{i_d})$ such that $V_{(n_1)} \preceq \cdots \preceq V_{(n_d)}$, according to an appropriate linear ordering $\preceq$ in $S_+^{d-1}$.

The proof that $\Psi_0(B_r^d \times B)$ is asymptotically normal as $r \to \infty$ relies on the following basic property of the Poisson process $\Psi$. Given the number $N_r = \Psi([-r, r] \times S_+^{d-1})$ of hyperplanes hitting $B_r^d$, say, $N_r = n$, the random vectors $X_i = (P_i, V_i), i = 1, \ldots, n$, are i.i.d. (and also independent of $N_r$) with independent components, where $P_i$ is uniformly distributed on $[-r, r]$ and $V_i$ has the distribution $\Theta$. Notice that $N_r$ is Poisson distributed with mean $2\lambda r$, which corresponds to (2.2) for $k = d - 1$.

In this way, we get that

$$\Psi_0(B_r^d \times B) \stackrel{\mathrm{d}}{=} \frac{1}{d!} \sum_{1 \leq i_1, \ldots, i_d \leq N_r}^* f_B(X_{i_1}, \ldots, X_{i_d}), \tag{3.1}$$

where $\stackrel{\mathrm{d}}{=}$ means equality in distribution, the sum $\sum^*$ runs over pairwise distinct indices, and

$$\begin{aligned}f_B((p_1, v_1), \ldots, (p_d, v_d)) \\ = \chi(H(p_1, v_1) \cap \cdots \cap H(p_d, v_d) \cap B_r^d) \mathbb{1}_B(v_{(1)}, \ldots, v_{(d)}),\end{aligned} \tag{3.2}$$

where $\chi(K) = 1$ for $K \neq \varnothing$ and $\chi(\varnothing) = 0$. Since the function $f_B : (\mathbb{R} \times S_+^{d-1})^d \to \{0, 1\}$ is symmetric and measurable, the right-hand side of (3.1) divided by $\binom{N_r}{d}$ and conditioned on $N_r = n$ is a $U$-statistic of order $d$ with kernel function $f = f_B$ as defined in (2.8).

3.1. *Moment formulae.* Since the first components $P_i$ of the i.i.d. random vectors $X_i = (P_i, V_i), i \geq 1$, are uniformly distributed on $[-r, r]$, the expectations $\mathbb{E} f_B(X_1, \ldots, X_d)$ do not depend on $r > 0$. To simplify notation, we put

$$\begin{aligned}\sigma_B^{(j,d)} &= \mathbb{E}(f_B(X_1, \ldots, X_d) f_B(X_{d-j+1}, \ldots, X_{2d-j})) \\ &= \mathbb{E} g_B^2(X_{d-j+1}, \ldots, X_d),\end{aligned} \tag{3.3}$$

where $g_B((p_1,v_1),\ldots,(p_j,v_j)) = \mathbb{E} f_B((p_1,v_1),\ldots,(p_j,v_j), X_{j+1},\ldots,X_d)$ for $j=1,\ldots,d$. Notice that also the second moments $\sigma_B^{(j,d)}$ do not depend on $r>0$. Now we formulate a first auxiliary result.

LEMMA 3.1. *For any* $B \in \mathcal{B}((S_+^{d-1})^d)$,

(3.4) $$\mathbb{E}\Psi_0(B_r^d \times B) = \frac{(2\lambda r)^d}{d!}\mathbb{E} f_B(X_1,\ldots,X_d)$$

*and*

(3.5) $$\operatorname{Var}\Psi_0(B_r^d \times B) = \sum_{j=1}^d \frac{(2\lambda r)^{2d-j}}{j!((d-j)!)^2}\sigma_B^{(j,d)}.$$

PROOF. By the symmetry of the function $f_B$ defined in (3.2) combined with the independence between $N_r$ and the i.i.d. sequence $X_i = (P_i, V_i), i \geq 1$, it is easily seen from (3.1) that

$$\mathbb{E}\Psi_0(B_r^d \times B) = \frac{1}{d!}\mathbb{E}(N_r(N_r-1)\cdots(N_r-d+1))\mathbb{E} f_B(X_1,\ldots,X_d).$$

Note that the $d$th factorial moment of a Poisson distributed random variable is equal to the $d$th power of its mean. Thus,

(3.6) $$\mathbb{E}(N_r(N_r-1)\cdots(N_r-d+1)) = (\mathbb{E} N_r)^d = (2\lambda r)^d$$

which proves (3.4). To derive a formula for the variance $\operatorname{Var}\Psi_0(B_r^d \times B)$, we again utilize the symmetry of $f_B$ and employ some simple combinatorial arguments which lead to

$$\mathbb{E}(\Psi_0(B_r^d \times B))^2$$
$$= \frac{1}{(d!)^2}\mathbb{E}\left(\sum_{1\leq i_1,\ldots,i_d \leq N_r}^* f_B(X_{i_1},\ldots,X_{i_d})\right)^2$$
$$= \sum_{j=0}^d \frac{j!}{(d!)^2}\binom{d}{j}\binom{d}{j}$$
$$\quad \times \mathbb{E}\left(\sum_{1\leq i_1,\ldots,i_{2d-j}\leq N_r}^* f_B(X_{i_1},\ldots,X_{i_d})f_B(X_{i_{d-j+1}},\ldots,X_{i_{2d-j}})\right)$$
$$= \sum_{j=0}^d \frac{\mathbb{E}(N_r(N_r-1)\cdots(N_r-2d+j+1))}{j!((d-j)!)^2}$$
$$\quad \times \mathbb{E}(f_B(X_1,\ldots,X_d)f_B(X_{d-j+1},\ldots,X_{2d-j})).$$

Finally, applying (3.6) with $d$ replaced by $2d-j$ for $j = 0, 1, \ldots, d$, and noting that the summand for $j = 0$ in the last line coincides with $(\mathbb{E}\Psi_0(B_r^d \times B))^2$, we obtain (3.5), which completes the proof of Lemma 3.1.  □

Notice that, as an immediate consequence of (3.5), we obtain the limiting relation

$$(3.7) \qquad \lim_{r \to \infty} \frac{\operatorname{Var} \Psi_0(B_r^d \times B)}{r^{2d-1}} = \frac{(2\lambda)^{2d-1}}{((d-1)!)^2} \sigma_B^{(1,d)}.$$

3.2. *Central limit theorem for the number of intersection points.* We now are in a position to formulate and prove a CLT for the number $\Psi_0(B_r^d \times B)$ of marked intersection points as $r \to \infty$, where the centering and normalizing constants have been derived in Lemma 3.1 and in (3.7), respectively.

THEOREM 3.1. *Let $B \in \mathcal{B}((S_+^{d-1})^d)$ be chosen such that $\sigma_B^{(1,d)} > 0$. Then,*

$$\frac{(d-1)!}{(2\lambda r)^{d-1/2}} \left( \Psi_0(B_r^d \times B) - \frac{(2\lambda r)^d}{d!} \mathbb{E} f_B(X_1, \ldots, X_d) \right) \xrightarrow[r \to \infty]{\mathrm{d}} \mathcal{N}(0, \sigma_B^{(1,d)}).$$

PROOF. Note that (3.1) is equivalent to the equality $\Psi_0(B_r^d \times B) \stackrel{\mathrm{d}}{=} \binom{N_r}{d} U_{N_r}^{(d)}(f_B)$. For any fixed $N_r = n \geq d$, the random multiple sum $U_{N_r}^{(d)}(f_B)$ is a $U$-statistic as defined in (2.8) of order $d$ with kernel function $f_B$ defined in (3.2). Let $n_r$ denote the expected value $\mathbb{E} N_r = 2\lambda r$ and let $\mu_B = \mathbb{E} f_B(X_1, \ldots, X_d)$. Then, Hoeffding's decomposition (2.9) yields

$$\Psi_0(B_r^d \times B) - \frac{(2\lambda r)^d}{d!} \mathbb{E} f_B(X_1, \ldots, X_d)$$

$$\stackrel{\mathrm{d}}{=} \left( \binom{N_r}{d} - \frac{n_r^d}{d!} \right) \mu_B + \binom{N_r}{d} \frac{d}{N_r} \sum_{i=1}^{N_r} (g_B(X_i) - \mu_B) + \binom{N_r}{d} R_{N_r}^{(d)}(f_B)$$

$$= \left( \binom{N_r}{d} - N_r \binom{N_r - 1}{d-1} + n_r \binom{N_r - 1}{d-1} - \frac{n_r^d}{d!} \right) \mu_B + \binom{N_r}{d} R_{N_r}^{(d)}(f_B)$$

$$+ \binom{N_r - 1}{d-1} \left( \sum_{i=1}^{N_r} g_B(X_i) - n_r \mu_B \right),$$

where $g_B(x) = \mathbb{E} f_B(x, X_2, \ldots, X_d)$ for $x \in [-r, r] \times S_+^{d-1}$. Since $N_r$ is independent of $X_1, X_2, \ldots$, the estimate (2.10) implies

$$\mathbb{E}((R_{N_r}^{(d)})^2 | N_r = n) \leq \frac{c_d}{n^2} \mathbb{E} f_B^2(X_1, \ldots, X_d) \qquad \text{for } n \geq d.$$

Hence, we conclude that

$$\mathbb{E}\left(\binom{N_r}{d} R_{N_r}^{(d)}(f_B)\right)^2 = \sum_{n \geq d} \binom{n}{d}^2 \mathbb{E}((R_{N_r}^{(d)})^2 | N_r = n) \mathbb{P}(N_r = n)$$

$$\leq \frac{c_d \mathbb{E} f_B^2(X_1, \ldots, X_d)}{d^2} \mathbb{E}\left(\binom{N_r - 1}{d - 1}\right)^2.$$

To determine the second moment of the binomial coefficient $\binom{N_r-1}{d-1}$, we use the expansions

$$x(x-1) \cdots (x-k+1) = \sum_{j=1}^{k} s_{j,k}^{(1)} x^j$$

and

$$x^k = \sum_{j=1}^{k} s_{j,k}^{(2)} x(x-1) \cdots (x-j+1),$$

where $s_{j,k}^{(1)}$ and $s_{j,k}^{(2)}$ denote the Stirling numbers of the first and second kind, respectively; see, for example, [26]. From (3.6), it is seen that $\mathbb{E} N_r^k$ is a polynomial of degree $k$ in $n_r$. Furthermore, $\mathbb{E}((N_r - 1)(N_r - 2) \cdots (N_r - d + 1))^2$ can be expressed as a polynomial of degree $2d - 2$ in $n_r$ such that

$$\frac{1}{r^{2d-1}} \mathbb{E}\left(\binom{N_r - 1}{d - 1}\right)^2 \xrightarrow[r \to \infty]{} 0.$$

Hence,

(3.8)  $$\frac{1}{r^{d-1/2}} \binom{N_r}{d} R_{N_r}^{(d)}(f_B) \xrightarrow[r \to \infty]{\mathbb{P}} 0,$$

where $\xrightarrow{\mathbb{P}}$ denotes convergence in probability. Next we show that

(3.9)  $$\frac{1}{r^{d-1/2}}\left(\binom{N_r}{d} - N_r \binom{N_r - 1}{d - 1} + n_r \binom{N_r - 1}{d - 1} - \frac{n_r^d}{d!}\right) \xrightarrow[r \to \infty]{\mathbb{P}} 0.$$

By virtue of $\frac{N_r}{r} \xrightarrow[r \to \infty]{\mathbb{P}} 2\lambda = \frac{n_r}{r}$ and since the difference $\binom{N_r}{d} - \frac{N_r^d}{d!}$ can be written as a polynomial of degree $d - 1$ in $N_r$ with coefficients not depending on $r$, it is easily seen that both

$$\frac{1}{r^{d-1/2}}\left(\binom{N_r}{d} - \frac{N_r^d}{d!}\right) \xrightarrow[r \to \infty]{\mathbb{P}} 0$$

and

$$\frac{N_r - n_r}{r^{d-1/2}}\left(\binom{N_r - 1}{d - 1} - \frac{N_r^{d-1}}{(d - 1)!}\right) \xrightarrow[r \to \infty]{\mathbb{P}} 0.$$

Similarly, after some elementary manipulations, we find that

$$\frac{1}{r^{d-1/2}}(N_r^d - n_r^d - d(N_r - n_r)N_r^{d-1})$$

$$= \frac{N_r - n_r}{r^{d-1/2}}\left(\sum_{k=0}^{d-1} n_r^k N_r^{d-1-k} - dN_r^{d-1}\right)$$

$$= -\left(\frac{N_r - n_r}{r^{3/4}}\right)^2 \frac{1}{r^{d-2}}\left(\sum_{k=0}^{d-2}(d-k-1)n_r^k N_r^{d-k-2}\right) \xrightarrow[r\to\infty]{\mathbb{P}} 0,$$

which in combination with the previous relation proves (3.9). Combining (3.8), (3.9) and

$$(3.10) \qquad \frac{1}{r^{d-1}}\binom{N_r - 1}{d - 1} \xrightarrow[r\to\infty]{\mathbb{P}} \frac{(2\lambda)^{d-1}}{(d-1)!},$$

and applying Slutsky's lemma (see, e.g., [14]), we see that the subsequent Lemma 3.2 completes the proof. □

LEMMA 3.2. *Under the conditions of Theorem* 3.1,

$$\frac{1}{\sqrt{2\lambda r}}\left(\sum_{i=1}^{N_r} g_B(X_i) - 2\lambda r \mathbb{E} f_B(X_1, \ldots, X_d)\right) \xrightarrow[r\to\infty]{d} \mathcal{N}(0, \sigma_B^{(1,d)}).$$

PROOF. Let again $n_r = \mathbb{E}N_r = 2\lambda r$ and let, furthermore, $\mu_B = \mathbb{E}g_B(X_1) = \mathbb{E}f_B(X_1, \ldots, X_d)$. The characteristic function of

$$\xi_r = n_r^{-1/2}\left(\sum_{i=1}^{N_r} g_B(X_i) - n_r \mu_B\right)$$

is then given by

$$\mathbb{E}e^{it\xi_r} = \exp(-it\sqrt{n_r}\mu_B)\mathbb{E}\exp\left(\frac{it}{\sqrt{n_r}}\sum_{i=1}^{N_r} g_B(X_i)\right).$$

The characteristic function on the right-hand side can be simplified by the fact that $N_r$ is independent of the sequence $X_1, X_2, \ldots$ and that the probability generating function $\mathbb{E}z^{N_r}$ takes the form $\exp(n_r(z-1))$ for any complex $z$. Thus,

$$\mathbb{E}e^{it\xi_r} = \exp\left(n_r \mathbb{E}\left(e^{(it/\sqrt{n_r})g_B(X_1)} - 1 - \frac{it}{\sqrt{n_r}}g_B(X_1)\right)\right)$$

or, equivalently, $\log \mathbb{E}e^{it\xi_r}$ is given by

$$n_r \mathbb{E}\left(\exp\left(\frac{it}{\sqrt{n_r}}g_B(X_1)\right) - 1 - \frac{it}{\sqrt{n_r}}g_B(X_1) - \frac{(it)^2}{2n_r}g_B^2(X_1)\right) - \frac{t^2 \sigma_B^{(1,d)}}{2}.$$

The well-known inequality $|e^{\mathrm{i}x} - 1 - \mathrm{i}x - \frac{(\mathrm{i}x)^2}{2}| \leq \frac{|x|^3}{6}$ for any $x \in \mathbb{R}$ combined with $n_r^{-1/2}\mathbb{E}|g_B^3(X_1)| \underset{r\to\infty}{\longrightarrow} 0$ gives

$$\log \mathbb{E} e^{\mathrm{i}t\xi_r} \underset{r\to\infty}{\longrightarrow} -\frac{t^2 \sigma_B^{(1,d)}}{2} \qquad \text{for } t \in \mathbb{R},$$

which is equivalent to the assertion of Lemma 3.2. □

**4. Extensions and applications of Theorem 3.1.** The results of Section 3 can be generalized as follows. Let again $\Psi = \sum_{i\geq 1} \delta_{(P_i, V_i)}$ be the marked-point-process representation of a stationary Poisson hyperplane process $\Phi$ in $\mathbb{R}^d$ with intensity $\lambda > 0$ and with spherical orientation distribution $\Theta$ such that $\Theta(L \cap S_+^{d-1}) < 1$ for any $L \in \mathcal{L}_{d-1}^d$. Notice that the expositions in this section are, for notational ease, presented for the particular case $B = (S_+^{d-1})^d$ only.

Instead of the intersection of $d$ hyperplanes, we consider, in Section 4.1, a generalized version of Theorem 3.1. In particular, for any fixed $k \in \{0, \ldots, d-1\}$, we define families of statistics $\Psi_k(B_r^d)$ and $\zeta_k(B_r^d)$ for $k = 0, \ldots, d-1$, denoting the number of $k$-flats hitting $B_r^d$ and their total $k$-volume in $B_r^d$, respectively. These statistics, closely related to the $k$-flat intersection process $\Phi_k$ defined in (2.5), have the form of $U$-statistics of order $d - k$ without normalizing factor $\binom{N_r}{d-k}$, that is,

$$\Psi_k(B_r^d) = \Phi_k(\{L \in \mathcal{A}_k^d : L \cap B_r^d \neq \varnothing\})$$
(4.1)
$$\stackrel{\mathrm{d}}{=} \frac{1}{(d-k)!} \sum_{1 \leq i_1, \ldots, i_{d-k} \leq N_r}^* \chi\left(\bigcap_{j=1}^{d-k} H(X_{i_j}) \cap B_r^d\right)$$

and

$$\zeta_k(B_r^d) = \sum_{L \in \mathrm{supp}(\Phi_k)} \nu_k(B_r^d \cap L)$$
(4.2)
$$\stackrel{\mathrm{d}}{=} \frac{1}{(d-k)!} \sum_{1 \leq i_1, \ldots, i_{d-k} \leq N_r}^* \nu_k\left(\bigcap_{j=1}^{d-k} H(X_{i_j}) \cap B_r^d\right).$$

In Section 4.3 we prove a multivariate CLT for $d$-dimensional vectors consisting of these, suitably normalized, random variables.

4.1. *CLTs for point processes of $k$-flats.* In analogy to Section 3.1, we first note that the expectations $\mathbb{E}\chi(H(X_1) \cap \cdots \cap H(X_{d-k}) \cap B_r^d)$ and $r^{-k}\mathbb{E}\nu_k(H(X_1) \cap \cdots \cap H(X_{d-k}) \cap B_r^d)$ do not depend on $r > 0$ since the first components $P_i$ of the i.i.d. random vectors $X_i = (P_i, V_i), i \geq 1$, are uniformly distributed on $[-r, r]$.

As an extension of Lemma 3.1, for $B = (S_+^{d-1})^d$, we now determine the first-order and second-order moments of the random variables $\Psi_k(B_r^d)$ and $\zeta_k(B_r^d)$ for $k = 0, \ldots, d-1$. For this reason, let

$$\sigma_{\chi,k}^{(1,d-k)} = \mathbb{E}\left(\chi\left(\bigcap_{i=1}^{d-k} H(X_i) \cap B_r^d\right)\chi\left(\bigcap_{i=d-k}^{2d-2k-1} H(X_i) \cap B_r^d\right)\right) = \mathbb{E}g_{\chi,k}^2(X_{d-k})$$

and

$$\sigma_{\nu,k}^{(1,d-k)} = r^{-2k}\mathbb{E}\left(\nu_k\left(\bigcap_{i=1}^{d-k} H(X_i) \cap B_r^d\right)\nu_k\left(\bigcap_{i=d-k}^{2d-2k-1} H(X_i) \cap B_r^d\right)\right)$$
$$= \mathbb{E}g_{\nu,k}^2(X_{d-k}),$$

where

$$(4.3) \quad g_{\chi,k}((p,v)) = \mathbb{E}\chi(H(X_1) \cap \cdots \cap H(X_{d-k-1}) \cap H(p,v) \cap B_r^d)$$

and

$$(4.4) \quad g_{\nu,k}((p,v)) = r^{-k}\mathbb{E}\nu_k(H(X_1) \cap \cdots \cap H(X_{d-k-1}) \cap H(p,v) \cap B_r^d)$$

for $(p,v) \in [-r,r] \times S_+^{d-1}$. Notice that the second moments $\sigma_{\chi,k}^{(1,d-k)}$ and $\sigma_{\nu,k}^{(1,d-k)}$ do also not depend on $r > 0$. Using this notation, we can state the following moment formulae.

LEMMA 4.1. *For each $k = 0, \ldots, d-1$,*

$$(4.5) \quad \mathbb{E}\Psi_k(B_r^d) = \frac{(2\lambda r)^{d-k}}{(d-k)!}\mathbb{P}(H(X_1) \cap \cdots \cap H(X_{d-k}) \cap B_r^d \neq \varnothing),$$

$$(4.6) \quad \mathbb{E}\zeta_k(B_r^d) = \frac{(2\lambda r)^{d-k}}{(d-k)!}\mathbb{E}\nu_k(H(X_1) \cap \cdots \cap H(X_{d-k}) \cap B_r^d),$$

*and*

$$(4.7) \quad \lim_{r \to \infty} \frac{\operatorname{Var}\Psi_k(B_r^d)}{r^{2d-2k-1}} = \frac{(2\lambda)^{2d-2k-1}}{((d-k-1)!)^2}\sigma_{\chi,k}^{(1,d-k)},$$

$$(4.8) \quad \lim_{r \to \infty} \frac{\operatorname{Var}\zeta_k(B_r^d)}{r^{2d-1}} = \frac{(2\lambda)^{2d-2k-1}}{((d-k-1)!)^2}\sigma_{\nu,k}^{(1,d-k)}.$$

PROOF. In analogy to the proof of Lemma 3.1, we get for $k \in \{0, \ldots, d-1\}$ that

$$\mathbb{E}\Psi_k(B_r^d) = \mathbb{E}\binom{N_r}{d-k}\mathbb{E}\chi(H(X_1) \cap \cdots \cap H(X_{d-k}) \cap B_r^d)$$

and
$$\mathbb{E}\zeta_k(B_r^d) = \mathbb{E}\binom{N_r}{d-k}\mathbb{E}\nu_k(H(X_1)\cap\cdots\cap H(X_{d-k})\cap B_r^d).$$

Applying (3.6), we obtain that
$$\mathbb{E}\binom{N_r}{d-k} = \frac{(2\lambda r)^{d-k}}{(d-k)!},$$

which gives both (4.5) and (4.6). Furthermore, again arguing along the lines of the proof of Lemma 3.1, we obtain

$$\operatorname{Var}\Psi_k(B_r^d) = \sum_{j=1}^{d-k} \frac{(2\lambda r)^{2d-2k-j}}{j!((d-k-j)!)^2}$$
$$\times \mathbb{E}\left(\chi\left(\bigcap_{p=1}^{d-k} H(X_p)\cap B_r^d\right)\chi\left(\bigcap_{q=d-k-j+1}^{2(d-k)-j} H(X_q)\cap B_r^d\right)\right)$$

and

$$\operatorname{Var}\zeta_k(B_r^d) = \sum_{j=1}^{d-k} \frac{(2\lambda r)^{2d-2k-j}}{j!((d-k-j)!)^2}$$
$$\times \mathbb{E}\left(\nu_k\left(\bigcap_{p=1}^{d-k} H(X_p)\cap B_r^d\right)\nu_k\left(\bigcap_{q=d-k-j+1}^{2(d-k)-j} H(X_q)\cap B_r^d\right)\right).$$

Hence, after dividing by $r^{2d-2k-1}$ and $r^{2d-1}$, respectively, and letting $r\to\infty$, we get the desired relationships (4.7) and (4.8). $\square$

Recall now that the random variables $\Psi_k(B_r^d)$ and $\zeta_k(B_r^d)$ given in (4.1) and (4.2), respectively, can be expressed as $U$-statistics, allowing for Hoeffding's decomposition (2.9) to be applied. Hence, we can state the following CLTs, the proofs of which are in complete analogy to the proof of Theorem 3.1 and are therefore omitted.

THEOREM 4.1. *Let $\Phi$ be a Poisson hyperplane process with intensity $\lambda > 0$ and nondegenerate spherical orientation distribution $\Theta$. Then, for $k = 0, 1, \ldots, d-1$,*

(4.9) $\quad Z_{k,r}^{(d)}(\chi) = \dfrac{(d-k-1)!}{(2\lambda r)^{d-k-1/2}}(\Psi_k(B_r^d) - \mathbb{E}\Psi_k(B_r^d)) \xrightarrow[r\to\infty]{\mathrm{d}} \mathcal{N}(0, \sigma_{\chi,k}^{(1,d-k)})$

*and*

(4.10) $\quad Z_{k,r}^{(d)}(\nu) = \dfrac{(d-k-1)!r^{-k}}{(2\lambda r)^{d-k-1/2}}(\zeta_k(B_r^d) - \mathbb{E}\zeta_k(B_r^d)) \xrightarrow[r\to\infty]{\mathrm{d}} \mathcal{N}(0, \sigma_{\nu,k}^{(1,d-k)}).$

If the stationary Poisson hyperplane process $\Phi$ is also isotropic, that is, $\Theta$ is the uniform distribution on $S_+^{d-1}$, the intensity $\lambda_k$ of the $k$-flat intersection process $\Phi_k$ induced by $\Phi$ is given by (2.7) for $k = 0, \ldots, d-1$. Using this in the following Lemma 4.2, we establish explicit expressions for the asymptotic variances occurring in the CLTs (4.9) and (4.10).

LEMMA 4.2. *Let $\Theta$ be the uniform distribution on $S_+^{d-1}$ and let $k \in \{0, 1, \ldots, d-1\}$. Then,*

$$(4.11) \qquad \mathbb{E}\Psi_k(B_r^d) = \lambda_k \kappa_{d-k} r^{d-k} \quad \text{and} \quad \mathbb{E}\zeta_k(B_r^d) = \lambda_k \kappa_d r^d,$$

*with $\lambda_k$ given in (2.7). Moreover,*

$$(4.12) \qquad \sigma_{\chi,k}^{(1,d-k)} = \frac{(\kappa_{d-k-1}(d-k-1)!)^2}{(2d-2k-1)!} \left(\frac{d!\kappa_d}{k!\kappa_k}\right)^2 \left(\frac{\kappa_{d-1}}{d\kappa_d}\right)^{2(d-k)}$$

*and*

$$(4.13) \qquad \sigma_{\nu,k}^{(1,d-k)} = \frac{(2^k \kappa_{d-1}(d-1)!)^2}{(2d-1)!} \left(\frac{d!\kappa_d}{k!\kappa_k}\right)^2 \left(\frac{\kappa_{d-1}}{d\kappa_d}\right)^{2(d-k)}.$$

PROOF. Both mean values in (4.11) are an immediate consequence of (2.2) and (2.3), respectively, where only the stationarity of the Poisson hyperplane process $\Phi$ is necessary. However, in case $\Phi$ is additionally isotropic, the intensities $\lambda_k$ can be explicitely determind by (2.7). To show (4.12), we use the relation $\sigma_{\chi,k}^{(1,d-k)} = \mathbb{E}g_{\chi,k}^2(X_{d-k})$, where the function $g_{\chi,k}((p,v))$ is defined in (4.3). Since the random vectors $X_i = (P_i, V_i), i = 1, \ldots, d-k$, are i.i.d. with independent components, where $P_i$ is uniformly distributed on $[-r, r]$ and $V_i$ has the uniform distribution $\Theta$ on $S_+^{d-1}$, the function $g_{\chi,k}((p,v))$ may be written in the form

$$g_{\chi,k}((p,v)) = \frac{1}{(2r)^{d-k-1}}$$
$$\times \int_{S_+^{d-1}} \int_{\mathbb{R}} \cdots \int_{S_+^{d-1}} \int_{\mathbb{R}} \chi\left(\bigcap_{j=1}^{d-k-1} H(p_j, v_j)\right.$$
$$\left. \cap H(p,v) \cap B_r^d\right) dp_1 \, \Theta(dv_1)$$
$$\cdots dp_{d-k-1} \, \Theta(dv_{d-k-1}).$$

A closed expression for $g_{\chi,k}((p,v))$ is obtained by an iterated application of Crofton's formula

$$(4.14) \quad \int_{S_+^{d-1}} \int_{\mathbb{R}} V_j(K \cap H(p,v)) \, dp \, \Theta(dv) = \frac{\kappa_{d-1}}{d\kappa_d} \frac{(j+1)\kappa_{j+1}}{\kappa_j} V_{j+1}(K),$$

which holds for $j = 0, 1, \ldots, d-1$ and any convex compact set $K \subset \mathbb{R}^d$; see, for example, Corollary 3.3.2 in [29]. Here, $V_i(\cdot)$ denotes the $i$th intrinsic volume for $i = 0, 1, \ldots, d$, where $V_0(K) = \chi(K)$ and $V_d(K) = \nu_d(K)$ for any convex compact $K \subset \mathbb{R}^d$. Applying (4.14) successively for $j = 0, \ldots, d-k-1$, we get

$$g_{\chi,k}((p,v)) = \frac{(d-k-1)!\kappa_{d-k-1}}{(2r)^{d-k-1}}\left(\frac{\kappa_{d-1}}{d\kappa_d}\right)^{d-k-1} V_{d-k-1}(H(p,v) \cap B_r^d).$$

Since $H(p,v) \cap B_r^d$ is a $(d-1)$-dimensional ball with radius $\sqrt{r^2 - p^2}$, the invariance and homogeneity properties of $V_{d-k-1}(\cdot)$ yield

$$V_{d-k-1}(H(p,v) \cap B_r^d) = (r^2 - p^2)^{(d-k-1)/2} V_{d-k-1}(B_1^{d-1})$$
$$= (r^2 - p^2)^{(d-k-1)/2} \binom{d-1}{k} \frac{\kappa_{d-1}}{\kappa_k};$$

see [29], page 79. Summarizing the above steps, we arrive at

$$(4.15) \qquad g_{\chi,k}((p,v)) = \frac{\kappa_{d-k-1}}{2^{d-k-1}} \frac{d!\kappa_d}{k!\kappa_k} \left(\frac{\kappa_{d-1}}{d\kappa_d}\right)^{d-k} \left(1 - \frac{p^2}{r^2}\right)^{(d-k-1)/2}$$

for $(p,v) \in [-r,r] \times S_+^{d-1}$ and

$$\sigma_{\chi,k}^{(1,d-k)} = \frac{1}{2r} \int_{-r}^{r} \int_{S_+^{d-1}} g_{\chi,k}^2((p,v)) \Theta(dv)\, dp$$

$$= \left(\frac{d!\kappa_d \kappa_{d-k-1}}{k!\kappa_k 2^{d-k-1}}\right)^2 \left(\frac{\kappa_{d-1}}{d\kappa_d}\right)^{2(d-k)} \int_0^1 (1-p^2)^{d-k-1}\, dp.$$

Finally, we get (4.12) by observing that

$$(4.16) \qquad \int_0^1 (1 - p^2)^s\, dp = \frac{(s!2^s)^2}{(2s+1)!}, \qquad s = 0, 1, \ldots.$$

To verify (4.13), we make use of a formula for the second moment of $\zeta_k(B_r^d)$ obtained in [17], page 164:

$$(4.17) \quad \mathbb{E}\zeta_k^2(B_r^d) = \sum_{j=0}^{d-k} \frac{d!(d-j)!r^{2d-j}}{j!(k!(d-k-j)!)^2} \frac{\kappa_{2d-j}\kappa_d\kappa_{d-j}^3}{\kappa_{2(d-j)}\kappa_k^2}\left(\lambda\frac{\kappa_{d-1}}{d\kappa_d}\right)^{2(d-k)-j}$$

for $k = 0, \ldots, d-1$. From the second formula in (4.11) and (2.7), it is seen that the summand for $j = 0$ equals $(\mathbb{E}\zeta_k(B_r^d))^2$ and, therefore, in accordance with (4.8),

$$\sigma_{\nu,k}^{(1,d-k)} = \frac{((d-k-1)!)^2}{(2\lambda)^{2d-2k-1}} \lim_{r \to \infty} \frac{\operatorname{Var}\zeta_k(B_r^d)}{r^{2d-1}}$$

$$= \frac{d!(d-1)!}{(k!)^2} \frac{\kappa_{2d-1}\kappa_d\kappa_{d-1}^3}{\kappa_{2d-2}\kappa_k^2}\left(\frac{\kappa_{d-1}}{2d\kappa_d}\right)^{2d-2k-1}.$$

Finally, we obtain (4.13) by taking into account the relation

$$\frac{\kappa_{2d-1}}{\kappa_{2d-2}} = \frac{\Gamma(d)\sqrt{\pi}}{\Gamma(d+1/2)} = \frac{2^{2d-1}((d-1)!)^2}{(2d-1)!},$$

which follows from Legendre's duplication formula

(4.18) $$2^{2s-1}\Gamma(s+\tfrac{1}{2})\Gamma(s) = \sqrt{\pi}\Gamma(2s)$$

applied to the integer $s = d$. This completes the proof of Lemma 4.2. $\square$

Notice that (4.5), (4.6) and (4.11) yield simple relationships between $\mathbb{P}(H(X_1) \cap \cdots \cap H(X_{d-k}) \cap B_r^d \neq \varnothing)$, $\mathbb{E}\nu_k(H(X_1) \cap \cdots \cap H(X_{d-k}) \cap B_r^d)$, and $\lambda_k$ for $k = 0, 1, \ldots, d-1$. Furthermore, in case $\Phi$ is additionally isotropic, one can use (2.7) to get

$$\mathbb{P}(H(X_1) \cap \cdots \cap H(X_{d-k}) \cap B_r^d \neq \varnothing) = \frac{2^{k-d}d!}{k!} \frac{\kappa_d \kappa_{d-k}}{\kappa_k} \left(\frac{\kappa_{d-1}}{d\kappa_d}\right)^{d-k},$$

$$\mathbb{E}\nu_k(H(X_1) \cap \cdots \cap H(X_{d-k}) \cap B_r^d) = \frac{2^{k-d}d!}{k!} \frac{\kappa_d^2}{\kappa_k} \left(\frac{\kappa_{d-1}}{d\kappa_d}\right)^{d-k} r^k.$$

Notice that these formulae comply with results in [17], pages 160 and 161. Also, if we replace in the proof of (4.12) the function $g_{\chi,k}(\cdot)$ by $g_{\nu,k}(\cdot)$ defined in (4.4) (with $V_k = \nu_k$ instead of $V_0 = \chi$), it easily seen that

(4.19) $$g_{\nu,k}((p,v)) = \frac{\kappa_{d-1}}{2^{d-k-1}} \frac{d!\kappa_d}{k!\kappa_k} \left(\frac{\kappa_{d-1}}{d\kappa_d}\right)^{d-k} \left(1 - \frac{p^2}{r^2}\right)^{(d-1)/2}$$

for $(p,v) \in [-r,r] \times S_+^{d-1}$, which confirms once more (4.13) without using Matheron's formula (4.17). On the other hand, regarding (4.17) for $k = 0$ as a sum of power functions in $\lambda r$, we are able to determine the pair correlation function $g_0(r)$ of the stationary and isotropic point process $\Psi_0$ as a polynomial of degree $d-1$ in $(\lambda r)^{-1}$; see also [12]. More precisely, putting $g_0(r) = 1 + \sum_{j=1}^{d-1} c_{dj} (\lambda r)^{-j}$ and utilizing the relationship

$$\operatorname{Var} \Psi_0(B_r^d) = d\kappa_d \lambda_0^2 \int_0^{2r} \nu_d(B_r^d \cap (B_r^d + (u,0,\ldots,0)))(g_0(u) - 1)u^{d-1}\,du$$
$$+ \lambda_0 \kappa_d r^d$$

(see [32], page 131 for details), we get $c_{dj} = \binom{d-1}{j}(\frac{\kappa_{d-j}}{\kappa_d})^2(\frac{d\kappa_d}{\kappa_{d-1}})^j$ by comparison of coefficients. Here, we used (2.7) for $k = 0$ together with

$$\int_0^2 \nu_d(B_1^d \cap (B_1^d + (u,0,\ldots,0)))u^{d-j-1}\,du = \frac{1}{d\kappa_d} \int_{B_1^d} \int_{B_1^d} \frac{dx\,dy}{\|x-y\|^j}$$
$$= \frac{2^{d-j+1}\kappa_{2d-j}}{(d-j)(d-j+1)\kappa_{d-j}}$$

for $j = 1, \ldots, d-1$; see [29], page 177.

4.2. *Asymptotic confidence intervals for the k-flat intensities.* Let $\Phi$ be a stationary Poisson hyperplane process with intensity $\lambda > 0$ and nondegenerate orientation distribution $\Theta$. From the view point of spatial statistics, $\widehat{\lambda}_{k,r} = \Psi_k(B_r^d)/\nu_{d-k}(B_r^{d-k})$, as well as $\widetilde{\lambda}_{k,r} = \zeta_k(B_r^d)/\nu_d(B_r^d)$, are unbiased estimators for the intensity $\lambda_k$ of the stationary $k$-flat intersection process generated by $\Phi$; see Section 2.1. We mention that both estimators are strongly consistent since $\Phi$ is ergodic and even mixing; see [30], Chapter 6.4.

If $\Phi$ is additionally isotropic, then we know from (2.7) that $\lambda_k = a_{d,k}\lambda^{d-k}$. Together with (4.7) and (4.12), as well as (4.8) and (4.13), we obtain

$$\lim_{r\to\infty} \frac{r\operatorname{Var}\widehat{\lambda}_{k,r}}{4(d-k)^2} = \lambda^{2d-2k-1}a_{d,k}^2 b_{d-k}^2 \quad \text{and} \quad \lim_{r\to\infty} \frac{r\operatorname{Var}\widetilde{\lambda}_{k,r}}{4(d-k)^2} = \lambda^{2d-2k-1}a_{d,k}^2 b_d^2$$

for $k = 0, \ldots, d-1$, where

$$a_{d,k} = \binom{d}{k}\frac{\kappa_d}{\kappa_k}\left(\frac{\kappa_{d-1}}{d\kappa_d}\right)^{d-k} \quad \text{and} \quad b_j = \frac{2^{j-1}(j-1)!\kappa_{j-1}}{\sqrt{2(2j-1)!\kappa_j}} \quad \text{for } j = 1, 2, \ldots.$$

By means of (4.18), one can verify the inequality $b_j < b_{j+1}$ for $j = 1, 2, \ldots$, which in turn implies that

$$\lim_{r\to\infty} r\operatorname{Var}\widehat{\lambda}_{k,r} < \lim_{r\to\infty} r\operatorname{Var}\widetilde{\lambda}_{k,r} \quad \text{for } k = 1, \ldots, d-1.$$

Therefore, we prefer the estimators $\widehat{\lambda}_{k,r}$ to construct confidence intervals for $\lambda_k$. Notice that efficiency and other optimality properties of intensity estimators for stationary $k$-flat processes observed in fixed convex sampling windows have been studied in [27].

By Theorem 4.1, the estimators $\widehat{\lambda}_{k,r}$ and $\widetilde{\lambda}_{k,r}$ are asymptotically normally distributed. For example, together with the above abbreviations, (4.9) can be formulated as follows:

$$\sqrt{r}(\widehat{\lambda}_{k,r} - \lambda_k) \xrightarrow[r\to\infty]{\mathrm{d}} \mathcal{N}(0, 4(d-k)^2\lambda_k^{2-1/(d-k)}a_{d,k}^{1/(d-k)}b_{d-k}^2).$$

Next, we apply a variance-stabilizing transformation $f(x)$ for $x \geq 0$ to the latter CLT such that $\sqrt{r}(f(\widehat{\lambda}_{k,r}) - f(\lambda_k))$ has a Gaussian limit with mean 0 and variance 1; see, for example, [3] for details. It is easily checked that $f(x) = (x/a_{d,k})^{1/2(d-k)}/b_{d-k}$ is a suitable choice for such transformation, which gives rise to a $100(1-\alpha)\%$ (asymptotically exact) confidence interval $I_{k,r}^{(d)}(\alpha)$ for $\lambda_k$.

COROLLARY 4.2. *Under the above assumptions,*

$$\frac{\sqrt{r}}{b_{d-k}}a_{d,k}^{-1/(2(d-k))}(\widehat{\lambda}_{k,r}^{1/(2(d-k))} - \lambda_k^{1/(2(d-k))}) \xrightarrow[r\to\infty]{\mathrm{d}} \mathcal{N}(0, 1),$$

which is equivalent to $\mathbb{P}(\lambda_k \in I_{k,r}^{(d)}(\alpha)) \xrightarrow[r \to \infty]{} 1 - \alpha$ for any $\alpha \in (0,1)$ and $k = 0, \ldots, d-1$, where the interval $I_{k,r}^{(d)}(\alpha)$ is given by

$$\left[\left((\widehat{\lambda}_{k,r})^{1/(2(d-k))} - \frac{a_{d,k}^{1/(2(d-k))} b_{d-k}}{\sqrt{r}} z_{1-\alpha/2}\right)^{2(d-k)},\right.$$

$$\left.\left((\widehat{\lambda}_{k,r})^{1/(2(d-k))} + \frac{a_{d,k}^{1/(2(d-k))} b_{d-k}}{\sqrt{r}} z_{1-\alpha/2}\right)^{2(d-k)}\right]$$

and $z_{1-\alpha/2}$ denotes the $(1-\alpha/2)$-quantile of the standard normal distribution. Moreover, relation $\lambda = (\lambda_k/a_{d,k})^{1/(d-k)}$ permits to transform the interval $I_{k,r}^{(d)}(\alpha)$ into the interval

$$J_{k,r}^{(d)}(\alpha) = \left[\left(\left(\frac{\widehat{\lambda}_{k,r}}{a_{d,k}}\right)^{1/(2(d-k))} - \frac{b_{d-k}}{\sqrt{r}} z_{1-\alpha/2}\right)^2,\right.$$

$$\left.\left(\left(\frac{\widehat{\lambda}_{k,r}}{a_{d,k}}\right)^{1/(2(d-k))} + \frac{b_{d-k}}{\sqrt{r}} z_{1-\alpha/2}\right)^2\right],$$

which covers the intensity $\lambda$ of the hyperplane process $\Phi$ with probability $1 - \alpha$ as $r \to \infty$.

The results of Corollary 4.2 can be applied, for example, in the following context of telecommunication network modeling.

Figure 2 shows the road system of Paris, where the real data, given by $(x,y)$-coordinates describing the location of points which are interpreted as intersections of roads, have been connected in order to form (nonconvex) cells. Beyond, each data point is equipped with a mark displaying the type of road to which it belongs, for example, main roads and side streets. The main road system is of prior interest in the civil engineering part of strategic network planning since these roads gather extensions of cables from different side streets. Therefore, expensive network components and technical devices are placed along main roads and capacity planning and cost analysis for the cables running along these types of roads is an important task.

The SSLM introduced in Section 1 allows us to model telecommunication networks. An important part of the SSLM is the geometry model, that is, the model that represents the infrastructure or road system of a certain (urban) environment. In [8], a fitting procedure, based on minimization of distance measures and Monte Carlo test techniques, is introduced, where tessellation models are fitted to real infrastructure data. It turns out that urban main road systems are often best represented by a Poisson line tessellation (Poisson line process).

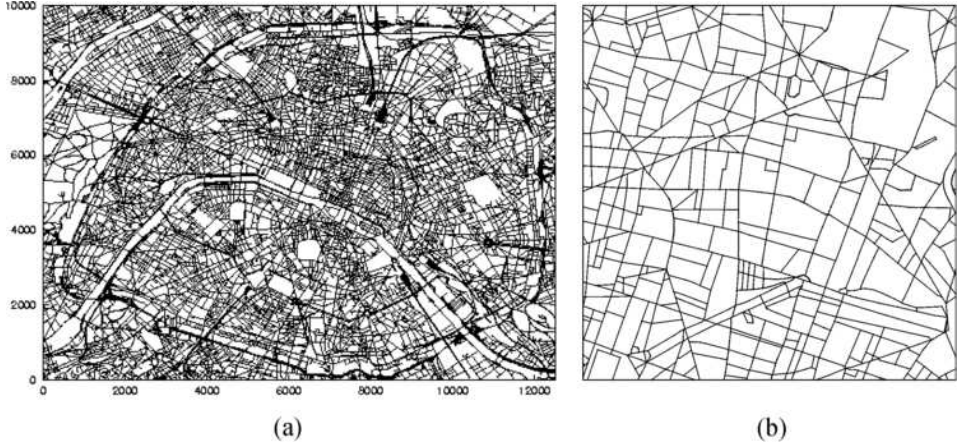

Fig. 2. *The road system of Paris:* (a) *city of Paris;* (b) *extracted data.*

From the civil engineering point of view, it is desirable to get information about the network's structure without having to apply too sophisticated models or methods. Furthermore, it is preferable for telecommunication engineers to get such information in the form of worst case and best case scenario values or to test their predictions. Our (asymptotic) confidence intervals, and, based upon them, our (asymptotic) tests, contribute to such a risk analysis of structural parameters of the network.

To take an example based on the observation of pointwise real data in a sufficiently large (spherical) sampling window $B_r^2$ and after application of the fitting procedure, it is possible to provide lower and upper bounds—$b_r^-(\alpha)$ and $b_r^+(\alpha)$, respectively—for the mean total length $\lambda$ of roads per unit area. To this end, consider the confidence interval $J_{0,r}^{(2)}(\alpha) = [b_r^-(\alpha), b_r^+(\alpha)]$ (from Corollary 4.2 for $d=2$ and $k=0$), where

$$b_r^\pm(\alpha) = \frac{1}{r}\bigg((\Psi_0(B_r^2))^{1/4} \pm \frac{2z_{1-\alpha/2}}{\pi\sqrt{3}}\bigg)^2.$$

This means that, based only on the knowledge of the number of road crossings within a large (spherical) region, $b_r^-(\alpha)\nu_2(W)$ and $b_r^+(\alpha)\nu_2(W)$ provide lower and upper bounds for the mean total length of the main road system intersecting a certain subregion $W$. These bounds can then be used to estimate lengths of cables and, in subsequent steps, to determine costs, as well as capacities of connection.

Moreover, based on the above confidence interval, one is able to test the null hypothesis $H_0: \lambda = \lambda^\star$, where $\lambda^\star$ denotes some specified value of the mean total length of lines per unit area, versus the alternative hypothesis $H_1: \lambda \neq \lambda^\star$. In the context of telecommunication, $\lambda^\star$ can be interpreted as the ratio of the total length of cables and the area (of an urban district) in which the cables are laid. This value for the cable length is often determined by practitioners according to their own rules of thumb. Based on the confidence interval $J_{0,r}^{(2)}(\alpha)$ ($r$ large enough), where $\alpha$ is a suitable significance level, $H_0$ would be rejected if

$$\widehat{\lambda}_{0,r} < \frac{1}{\pi}\left(\sqrt{\lambda^\star} - z_{1-\alpha/2}\frac{2}{\pi\sqrt{3r}}\right)^4 \quad \text{or} \quad \widehat{\lambda}_{0,r} > \frac{1}{\pi}\left(\sqrt{\lambda^\star} + z_{1-\alpha/2}\frac{2}{\pi\sqrt{3r}}\right)^4.$$

4.3. *Multivariate CLTs.* In this section we extend the results of Section 4.1 by establishing multivariate CLTs which describe the joint asymptotic behavior (as $r \to \infty$) of the closely correlated random variables $Z_{0,r}^{(d)}(\chi)$, ..., $Z_{d-1,r}^{(d)}(\chi)$, as well as of $Z_{0,r}^{(d)}(\nu), \ldots, Z_{d-1,r}^{(d)}(\nu)$, defined in (4.9) and (4.10), respectively. To begin with, we define the mixed second moments

$$\sigma_{\chi,k,l}^{(1,d-k,d-l)} = \mathbb{E}\left(\chi\left(\bigcap_{i=1}^{d-k} H(X_i) \cap B_r^d\right)\chi\left(\bigcap_{i=d-k}^{2d-k-l-1} H(X_i) \cap B_r^d\right)\right)$$

$$= \mathbb{E}(g_{\chi,k}(X_{d-k})g_{\chi,l}(X_{d-k}))$$

and

$$\sigma_{\nu,k,l}^{(1,d-k,d-l)} = r^{-(k+l)}\mathbb{E}\left(\nu_k\left(\bigcap_{i=1}^{d-k} H(X_i) \cap B_r^d\right)\nu_l\left(\bigcap_{i=d-k}^{2d-k-l-1} H(X_i) \cap B_r^d\right)\right)$$

$$= \mathbb{E}(g_{\nu,k}(X_{d-k})g_{\nu,l}(X_{d-k}))$$

for $k, l = 0, \ldots, d-1$, where the functions $g_{\chi,k}(\cdot)$ and $g_{\nu,l}(\cdot)$ are defined by (4.3) and (4.4), respectively. The following result generalizes Lemma 4.1.

LEMMA 4.3. *For any* $k, l = 0, \ldots, d-1$,

$$(4.20) \qquad \lim_{r\to\infty} \mathbb{E}(Z_{k,r}^{(d)}(\chi)Z_{l,r}^{(d)}(\chi)) = \sigma_{\chi,k,l}^{(1,d-k,d-l)}$$

*and*

$$(4.21) \qquad \lim_{r\to\infty} \mathbb{E}(Z_{k,r}^{(d)}(\nu)Z_{l,r}^{(d)}(\nu)) = \sigma_{\nu,k,l}^{(1,d-k,d-l)}.$$

PROOF. Both relations (4.20) and (4.21) can be shown by similar combinatorial arguments as employed in the proof of Lemma 3.1 to verify (3.5). Using the definition of $Z_{k,r}^{(d)}(\chi)$ as given in (4.9), we easily get that

$$\mathbb{E}(Z_{k,r}^{(d)}(\chi)Z_{l,r}^{(d)}(\chi)) = \frac{(d-k-1)!(d-l-1)!}{(2\lambda r)^{2d-k-l-1}}$$
$$\times (\mathbb{E}(\Psi_k(B_r^d)\Psi_l(B_r^d)) - \mathbb{E}\Psi_k(B_r^d)\mathbb{E}\Psi_l(B_r^d)).$$

According to (4.1), we may write for $0 \leq k \leq l \leq d-1$ that

$$\mathbb{E}(\Psi_k(B_r^d)\Psi_l(B_r^d))$$
$$= \sum_{j=0}^{d-l} \frac{\mathbb{E}(N_r(N_r-1)\cdots(N_r-(2d-k-l-j)+1))}{j!(d-k-j)!(d-l-j)!}$$
$$\times \mathbb{E}\left(\chi\left(\bigcap_{p=1}^{d-k} H(X_p) \cap B_r^d\right)\chi\left(\bigcap_{q=d-k-j+1}^{2d-k-l-j} H(X_q) \cap B_r^d\right)\right).$$

Using (3.6) and the fact that the summand for $j=0$ equals $\mathbb{E}\Psi_k(B_r^d)\mathbb{E}\Psi_l(B_r^d)$, we arrive at

$$\mathbb{E}(Z_{k,r}^{(d)}(\chi)Z_{l,r}^{(d)}(\chi))$$
$$= \frac{(d-k-1)!(d-l-1)!}{(2\lambda r)^{2d-k-l-1}}$$
$$\times \sum_{j=1}^{d-l} \frac{(2\lambda r)^{2d-k-l-j}}{j!(d-k-j)!(d-l-j)!}$$
$$\times \mathbb{E}\left(\chi\left(\bigcap_{p=1}^{d-k} H(X_p) \cap B_r^d\right)\chi\left(\bigcap_{q=d-k-j+1}^{2d-k-l-j} H(X_q) \cap B_r^d\right)\right).$$

Hence, dividing by $r^{2d-k-l-1}$ and letting $r \to \infty$ immediately yields (4.20). The proof of (4.21) is completely analogous. □

THEOREM 4.3. *Let the assumptions of Theorem 4.1 be satisfied. Then,*

(4.22) $$(Z_{k,r}^{(d)}(\chi))_{k=0}^{d-1} \xrightarrow[r\to\infty]{d} \mathcal{N}(o, \Sigma(\chi))$$

*and*

(4.23) $$(Z_{k,r}^{(d)}(\nu))_{k=0}^{d-1} \xrightarrow[r\to\infty]{d} \mathcal{N}(o, \Sigma(\nu)),$$

*where $\mathcal{N}(o, \Sigma(\chi))$ and $\mathcal{N}(o, \Sigma(\nu))$ are $d$-dimensional Gaussian mean $o = (0,\ldots,0)^\top$ random vectors with covariance matrices $\Sigma(\chi) = (\sigma_{\chi,k,l}^{(1,d-k,d-l)})_{k,l=0}^{d-1}$ and $\Sigma(\nu) = (\sigma_{\nu,k,l}^{(1,d-k,d-l)})_{k,l=0}^{d-1}$ with entries given by the limits (4.20) and (4.21), respectively.*

PROOF. Recall that due to the well-known Cramér–Wold device, the multivariate CLT (4.22) is equivalent to the one-dimensional CLT

(4.24) $$\sum_{k=0}^{d-1} t_k Z_{k,r}^{(d)}(\chi) \xrightarrow[r\to\infty]{d} \mathcal{N}(0, t^\top \Sigma(\chi) t),$$

for any $t = (t_0, \ldots, t_{d-1})^\top \in \mathbb{R}^d \setminus \{o\}$. This means that the proof of (4.22) can be put down to the case of the (one-dimensional) CLTs considered in Theorems 3.1 and 4.1. First, using (2.8), (4.1) and (4.9), we may rewrite the linear combination $\sum_{k=0}^{d-1} t_k Z_{k,r}^{(d)}(\chi)$ as follows:

$$\sum_{k=0}^{d-1} t_k Z_{k,r}^{(d)}(\chi) \stackrel{\mathrm{d}}{=} \sum_{k=0}^{d-1} t_k \frac{(d-k-1)!}{(2\lambda r)^{d-k-1/2}}$$
$$\times \left( U_{N_r}^{(d-k)}(\chi) \binom{N_r}{d-k} - \frac{(2\lambda r)^{d-k}}{(d-k)!} \mathbb{E} g_{\chi,k}(X_1) \right).$$

Next, we apply Hoeffding's decomposition (2.9) to the random $U$-statistic $U_{N_r}^{(d-k)}(\chi)$ with kernel function $\chi(\bigcap_{j=1}^{d-k} H(X_j) \cap B_r^d)$ and proceed in the same manner as in the proof of Theorem 3.1. In view of the limiting relations (3.8) and (3.9) with $d$ replaced by $d - k$ for $k = 0, 1, \ldots, d - 1$ and combined with Slutsky's lemma, we recognize that the weak limit of $\sum_{k=0}^{d-1} t_k Z_{k,r}^{(d)}(\chi)$ coincides with that of

$$\sum_{k=0}^{d-1} t_k \frac{(d-k-1)!}{(2\lambda r)^{d-k-1/2}} \binom{N_r - 1}{d-k-1} \left( \sum_{i=1}^{N_r} g_{\chi,k}(X_i) - 2\lambda r \mathbb{E} g_{\chi,k}(X_1) \right)$$

as $r \to \infty$. Finally, by means of (3.10) with $d$ again replaced by $d - k$ for $k = 0, 1, \ldots, d - 1$, it remains to show that

$$\frac{1}{\sqrt{2\lambda r}} \sum_{k=0}^{d-1} t_k \left( \sum_{i=1}^{N_r} g_{\chi,k}(X_i) - 2\lambda r \mathbb{E} g_{\chi,k}(X_1) \right)$$
$$\xrightarrow[r \to \infty]{\mathrm{d}} \mathcal{N}\left( 0, \sum_{k,l=0}^{d-1} t_k t_l \mathbb{E}(g_{\chi,k}(X_1) g_{\chi,l}(X_1)) \right).$$

However, the latter CLT is obtained by proving Lemma 3.2 once more for the function $\sum_{k=0}^{d-1} t_k g_{\chi,k}(\cdot)$ instead of $g_B(\cdot)$. To show the second assertion (4.23), we only need to repeat the just finished proof with the kernel function $r^{-k} \nu_k(\bigcap_{j=1}^{d-k} H(X_j) \cap B_r^d)$ and $g_{\nu,l}(\cdot)$ instead of $g_{\chi,l}(\cdot)$. This completes the proof of Theorem 4.3. □

As in Lemma 4.2, the additional assumption of isotropy allows to compute explicit formulae for the mixed second-order moments $\sigma_{\chi,k,l}^{(1,d-k,d-l)}$ and $\sigma_{\nu,k,l}^{(1,d-k,d-l)}$.

LEMMA 4.4. *Let $\Theta$ be the uniform distribution on $S_+^{d-1}$ and let $k, l \in \{0, \ldots, d-1\}$. Then,*

$$\sigma_{\chi,k,l}^{(1,d-k,d-l)} = \frac{(d!\kappa_d)^2 \kappa_{d-k-1} \kappa_{d-l-1}}{k! l! \kappa_k \kappa_l 2^{2d-k-l-1}} \frac{\kappa_{2d-k-l-1}}{\kappa_{2d-k-l-2}} \left( \frac{\kappa_{d-1}}{d\kappa_d} \right)^{2d-k-l}$$

(4.25)
$$= \sqrt{\sigma_{\chi,k}^{(1,d-k)}\sigma_{\chi,l}^{(1,d-l)}}\frac{\mathrm{B}((2d-k-l)/2,(2d-k-l)/2)}{\sqrt{\mathrm{B}(d-k,d-k)\mathrm{B}(d-l,d-l)}},$$

where $\mathrm{B}(s,t) = \int_0^1 x^{s-1}(1-x)^{t-1}\,dx = \Gamma(s)\Gamma(t)/\Gamma(s+t)$ denotes Euler's Beta function, and

(4.26)
$$\sigma_{\nu,k,l}^{(1,d-k,d-l)} = \frac{(\kappa_d\kappa_{d-1}d!(d-1)!)^2 2^{k+l}}{k!l!\kappa_k\kappa_l(2d-1)!}\left(\frac{\kappa_{d-1}}{d\kappa_d}\right)^{2d-k-l}$$
$$= \sqrt{\sigma_{\nu,k}^{(1,d-k)}\sigma_{\nu,l}^{(1,d-l)}}.$$

PROOF. Both (4.25) and (4.26) can be obtained using the shape of the functions $g_{\chi,k}(\cdot)$ and $g_{\nu,k}(\cdot)$ derived in the proof of Lemma 4.2. By (4.15) and the above definition of $\sigma_{\chi,k,l}^{(1,d-k,d-l)}$, we get that

$$\sigma_{\chi,k,l}^{(1,d-k,d-l)} = \frac{1}{2r}\int_{-r}^{r}\int_{S_+^{d-1}} g_{\chi,k}((p,v))g_{\chi,l}((p,v))\Theta(dv)\,dp$$
$$= \frac{d!\kappa_d\kappa_{d-k-1}}{k!\kappa_k 2^{d-k-1}}\frac{d!\kappa_d\kappa_{d-l-1}}{l!\kappa_l 2^{d-l-1}}\left(\frac{\kappa_{d-1}}{d\kappa_d}\right)^{2d-k-l}$$
$$\times \int_0^1(1-p^2)^{(2d-k-l-2)/2}\,dp.$$

Thus, by noting that

$$\int_0^1(1-p^2)^{(2d-k-l-2)/2}\,dp = \frac{\kappa_{2d-k-l-1}}{2\kappa_{2d-k-l-2}}$$
$$= 2^{2d-k-l-2}\mathrm{B}\left(\frac{2d-k-l}{2},\frac{2d-k-l}{2}\right)$$

(see, e.g., [29], page 80), we obtain the first part of (4.25), where the second identity in the previous line turns out to be a simple consequence of (4.18) for $s = (2d-k-l)/2$ and the very definition of the Beta function. The second part of (4.25) is seen by inserting the variances $\sigma_{\chi,k}^{(1,d-k)}$ given by (4.12) combined with $\mathrm{B}(d-k,d-k) = ((d-k-1)!)^2/(2d-2k-1)!$ for $k = 0, 1, \ldots, d-1$. Likewise, using (4.19), we get that

$$\sigma_{\nu,k,l}^{(1,d-k,d-l)} = \frac{1}{2r}\int_{-r}^{r}\int_{S_+^{d-1}} g_{\nu,k}((p,v))g_{\nu,l}((p,v))\Theta(dv)\,dp$$
$$= \frac{d!\kappa_d\kappa_{d-1}}{k!\kappa_k 2^{d-k-1}}\frac{d!\kappa_d\kappa_{d-1}}{l!\kappa_l 2^{d-l-1}}\left(\frac{\kappa_{d-1}}{d\kappa_d}\right)^{2d-k-l}\int_0^1(1-p^2)^{d-1}\,dp.$$

Hence, taking (4.16) for $s = d-1$, the first part of (4.26) is shown and the second equality is immediately seen from Lemma 4.2. $\square$

COROLLARY 4.4. *The covariance matrix $\Sigma(\chi)$ possesses always full rank $d$, whereas the rank of the covariance matrix $\Sigma(\nu)$ equals $1$ for any dimension $d \geq 1$. Moreover,*

$$(4.27) \qquad \frac{Z_{l,r}^{(d)}(\nu)}{\sqrt{\sigma_{\nu,l}^{(1,d-l)}}} - \frac{Z_{k,r}^{(d)}(\nu)}{\sqrt{\sigma_{\nu,k}^{(1,d-k)}}} \xrightarrow[r \to \infty]{\mathbb{P}} 0 \qquad \text{for } 0 \leq k < l \leq d-1.$$

PROOF. Notice that $\Sigma(\chi)$ possesses full rank if this matrix is strictly positive. This, however, can be seen since

$$\sum_{k,l=0}^{d-1} t_k t_l \mathrm{B}\left(\frac{2d-k-l}{2}, \frac{2d-k-l}{2}\right)$$

$$= \int_0^1 \left(\sum_{k=0}^{d-1} t_k (x(1-x))^{(d-k-1)/2}\right)^2 dx > 0$$

for any $(t_0, \ldots, t_{d-1})^\top \in \mathbb{R}^d \setminus \{o\}$, which means that the symmetric matrix with entries $\mathrm{B}((2d-k-l)/2, (2d-k-l)/2)$ for $k, l = 0, \ldots, d-1$ is strictly positive. In order to show that the rank of $\Sigma(\nu)$ is 1 for any $d \geq 1$, we only need to observe that the second equality in (4.26) implies that each entry of the asymptotic covariance matrix of the normalized random vector $(Z_{k,r}^{(d)}(\nu)/\sqrt{\sigma_{\nu,k}^{(1,d-k)}})_{k=0}^{d-1}$ equals 1. Then, (4.27) is a consequence of the structure of $\Sigma(\nu)$. □

The somewhat surprising result (4.27) states that the variance of the difference of any two components of $(Z_{k,r}^{(d)}(\nu)/\sqrt{\sigma_{\nu,k}^{(1,d-k)}})_{k=0}^{d-1}$ tends to zero as $r \to \infty$. Together with Slutsky's lemma this allows for the conclusion that the normal convergence in (4.10) for a single component, $Z_{0,r}^{(d)}(\nu)$, say, implies asymptotic normality of the other components. Thus, relation (4.27) can be interpreted as a kind of asymptotic second-order relationship for $k$-flat processes induced by stationary and isotropic Poisson hyperplane process. It would be of interest to see whether there is a pure geometric reasoning for (4.27).

**5. A special review of the planar case.** Throughout this section we assume that $d = 2$. Assuming isotropy of the underlying stationary Poisson line process, we first present a short review of Theorem 3.1 for the case where the ordered angles of the orientation vectors of intersecting pairs of lines are situated within a certain rectangle. In a second part of the present section we look at another type of a CLT for Poisson line processes, proven by Paroux [24], where the normalization is random. Applying directly Hoeffding's CLT (2.11) for $U$-statistics, we provide a new proof of Paroux's CLT

CLTS FOR POISSON HYPERPLANE TESSELLATIONS 27with random normalization, which has been derived in [24] by the "method of moments."

5.1. *Planar moment formulae.* We consider the marked-point-process representation $\Psi = \sum_{i \geq 1} \delta_{(P_i, V_i)}$ of a planar stationary and isotropic Poisson line process $\Phi$ with intensity $\lambda$. In this special case each orientation vector $V_i \in S^1_+$ is completely determined by the angle $\Gamma_i$ between the unit vector $V_i$ and the $x$-axis measured in anti-clockwise direction. Owing to isotropy, the angles $\Gamma_1, \Gamma_2, \ldots$ are independent and uniformly distributed on $[0, \pi]$. Therefore, $\mathrm{supp}(\Phi)$ consists of parameterized lines $\ell_{(P_i, \Gamma_i)}$ in $\mathbb{R}^2$ defined by $\ell_{(P_i, \Gamma_i)} = \{(x, y) \in \mathbb{R}^2 : x \cos \Gamma_i + y \sin \Gamma_i = P_i\}$.

For $r > 0$ fixed, $\Psi_0(B_r^2 \times B(a, b))$ is the random number of those intersection points $\ell_{(P_{n_1}, \Gamma_{n_1})} \cap \ell_{(P_{n_2}, \Gamma_{n_2})}$ in $B_r^2$ for which $(\Gamma_{(n_1)}, \Gamma_{(n_2)}) \in B(a, b)$, where $B(a, b) = [0, a] \times [0, b]$, $a \leq b$, is a rectangular subset of $[0, \pi]^2$. Recall that $\Psi_0(B_r^2 \times B(a, b))$ has the same distribution as the random double sum

$$\text{(5.1)} \qquad \sum_{1 \leq i < j \leq N_r} f_{B(a,b)}((P_i, \Gamma_i), (P_j, \Gamma_j)),$$

where, as in Section 3, the random vectors $(P_1, \Gamma_1), (P_2, \Gamma_2), \ldots : \Omega \to [-r, r] \times [0, \pi]$ are independent and uniformly distributed on $[-r, r] \times [0, \pi]$ and where $N_r : \Omega \to \{0, 1, \ldots\}$ is a Poisson distributed random variable with expectation $2\lambda r$ independent of the $(P_i, \Gamma_i)$'s. The function $f_{B(a,b)}$ given by

$$f_{B(a,b)}((p_1, \gamma_1), (p_2, \gamma_2)) = \chi(\ell_{(p_1, \gamma_1)} \cap \ell_{(p_2, \gamma_2)} \cap B_r^2) \mathbb{1}_{B(a,b)}(\gamma_{(1)}, \gamma_{(2)})$$

is symmetric since $(\gamma_{(1)}, \gamma_{(2)})$ are lexicographically ordered, that is, $(\gamma_{(1)}, \gamma_{(2)}) = (\gamma_1 \wedge \gamma_2, \gamma_1 \vee \gamma_2)$. With the abbreviations

$$\text{(5.2)} \qquad \begin{aligned} \mu_{B(a,b)} &= \mathbb{E} f_{B(a,b)}((P_1, \Gamma_1), (P_2, \Gamma_2)) \\ &= \mathbb{P}(\ell_{(P_1, \Gamma_1)} \cap \ell_{(P_2, \Gamma_2)} \cap B_1^2 \neq \varnothing, 0 \leq \Gamma_1 \wedge \Gamma_2 \leq a, \Gamma_1 \vee \Gamma_2 \leq b) \end{aligned}$$

and $\sigma^{(1,2)}_{B(a,b)}$ defined by (3.3), Theorem 3.1 claims that

$$\text{(5.3)} \quad \frac{1}{(2\lambda r)^{3/2}} \left( \Psi_0(B_r^2 \times B(a, b)) - \frac{(2r\lambda)^2}{2} \mu_{B(a,b)} \right) \xrightarrow[r \to \infty]{\mathrm{d}} \mathcal{N}(0, \sigma^{(1,2)}_{B(a,b)}),$$

where

$$\text{(5.4)} \qquad \begin{aligned} \mu_{B(a,b)} &= \frac{1}{4\pi} \int_0^a \int_0^b |\sin(u - v)| \, du \, dv \\ &= \frac{1}{4\pi} (2a - \sin a - \sin b + \sin(b - a)) \end{aligned}$$

and

$$\sigma_{B(a,b)}^{(1,2)} = \frac{2}{3\pi^3}(5a + 7(\sin(b-a) - \sin b)$$

(5.5)
$$+ (b + \sin b)\cos b$$

$$- (b - a + \sin(b-a))\cos(b-a)).$$

To verify (5.4), one can apply the general mean value formula (2.6) together with (2.7) for $d=2$ and $k=0$. However, we use (3.4) and a more direct approach to calculate the probability $\mu_{B(a,b)}$. By definition (5.2), we obtain

$$\mu_{B(a,b)} = \frac{1}{4\pi^2} \int_0^\pi \int_0^\pi \int_{-1}^1 \int_{-1}^1 \mathbb{1}_{\{x^2+y^2 \leq 1\}}(p_1, \gamma_1, p_2, \gamma_2)$$
$$\times \mathbb{1}_{[0,a]}(\gamma_1 \wedge \gamma_2)$$
$$\times \mathbb{1}_{[0,b]}(\gamma_1 \vee \gamma_2)\, dp_1\, dp_2\, d\gamma_1\, d\gamma_2,$$

where $(x,y)$ denotes the intersection point of the lines $\ell_{(p_1,\gamma_1)}$ and $\ell_{(p_2,\gamma_2)}$, that is,

$$x = -\frac{-p_1 \sin\gamma_2 + p_2 \sin\gamma_1}{\sin(\gamma_2 - \gamma_1)},$$
$$y = \frac{-p_1 \cos\gamma_2 + p_2 \cos\gamma_1}{\sin(\gamma_2 - \gamma_1)}.$$

Thus, by some elementary manipulations with trigonometric functions, we arrive at

$$\mu_{B(a,b)} = \frac{1}{2\pi^2} \int_0^a \int_0^b \int_{-1}^1 \sqrt{1-p_2^2}|\sin(\gamma_2 - \gamma_1)|\, dp_2\, d\gamma_1\, d\gamma_2,$$

which combined with $\int_{-1}^1 \sqrt{1-p^2}\, dp = \frac{\pi}{2}$ confirms (5.4). To determine the asymptotic variance $\sigma_{B(a,b)}^{(1,2)}$ appearing in (5.3), we argue similarly. By (3.3),

$$\sigma_{B(a,b)}^{(1,2)} = \mathbb{E}(f_{B(a,b)}((P_1,\Gamma_1),(P_2,\Gamma_2))f_{B(a,b)}((P_2,\Gamma_2),(P_3,\Gamma_3)))$$
$$= \frac{1}{8\pi^3}\int_0^\pi \int_0^\pi \int_0^\pi \int_{-1}^1 \int_{-1}^1 \int_{-1}^1 \mathbb{1}_{\{x^2+y^2 \leq 1\}}(p_1,\gamma_1,p_2,\gamma_2)$$
$$\times \mathbb{1}_{\{x^2+y^2 \leq 1\}}(p_2,\gamma_2,p_3,\gamma_3)$$
$$\times \mathbb{1}_{[0,a]}(\gamma_1 \wedge \gamma_2 \wedge \gamma_3)$$
$$\times \mathbb{1}_{[0,b]}(\gamma_1 \vee \gamma_2 \vee \gamma_3)\, dp_1\, dp_2\, dp_3\, d\gamma_1\, d\gamma_2\, d\gamma_3$$

$$= \frac{4}{8\pi^3} \int_0^b \int_0^b \int_0^b \int_{-1}^1 (1-p_2^2)|\sin(\gamma_2-\gamma_1)||\sin(\gamma_2-\gamma_3)|\,dp_2$$
$$\times \mathbb{1}_{[0,a]}(\gamma_1 \wedge \gamma_2 \wedge \gamma_3)\,d\gamma_1\,d\gamma_2\,d\gamma_3.$$

Using that $\int_{-1}^1 (1-p^2)\,dp = \frac{4}{3}$, a somewhat lengthy computation of the threefold integrals with respect to $\gamma_1, \gamma_2, \gamma_3$ in the latter line leads to formula (5.5).

5.2. *A CLT with random normalization.* The two-dimensional version (5.3) of our CLT for intersection points of Poisson line processes in $\mathbb{R}^2$ has close connections to a CLT with random normalization; see Theorem 3.1.3 in [24]. The proof given in [24] is based on the well-known "method of moments." The following Theorem 5.1 states the assertion of this CLT. However, using Hoeffding's CLT (2.11) for $U$-statistics, we obtain a much shorter proof.

THEOREM 5.1. *For arbitrary fixed $a, b \in [0, \pi]$ with $a < b$, let*

$$(5.6) \quad Z_{B(a,b)}^{(r)} = \frac{1}{(N_r(N_r-1))^{3/4}} \sum_{1 \leq i < j \leq N_r} (f_{B(a,b)}(X_i, X_j) - \mu_{B(a,b)}).$$

*Then,*

$$Z_{B(a,b)}^{(r)} \xrightarrow[r\to\infty]{d} \mathcal{N}(0, \sigma_{B(a,b)}^{(1,2)} - \mu_{B(a,b)}^2),$$

*where $\mu_{B(a,b)}$ and $\sigma_{B(a,b)}^{(1,2)}$ are given by (5.4) and (5.5).*

PROOF. To begin with, we rewrite $Z_{B(a,b)}^{(r)}$ as

$$Z_{B(a,b)}^{(r)} = \frac{N_r^{1/2}}{2}\left(1-\frac{1}{N_r}\right)^{1/4}\left(\frac{2}{N_r(N_r-1)}\sum_{1 \leq i < j \leq N_r}(f_{B(a,b)}(X_i,X_j)-\mu_{B(a,b)})\right),$$

where $f_{B(a,b)}(x_1, x_2)$ is a measurable, symmetric function in $x_1, x_2 \in [-r, r] \times S_+^1$ satisfying $\mathbb{E} f_{B(a,b)}(X_1, X_2) = \mu_{B(a,b)}$. Applying Hoeffding's CLT (2.11) for $U$-statistics yields

$$(5.7) \quad \frac{[n_r]^{1/2}}{2} U_{[n_r]}^{(2)}(f_{B(a,b)}) \xrightarrow[r\to\infty]{d} \mathcal{N}(0, \operatorname{Var} g_{B(a,b)}(X_1)),$$

where the $U$-statistic $U_{[n_r]}^{(2)}(f_{B(a,b)})$ is given by

$$U_{[n_r]}^{(2)}(f_{B(a,b)}) = \frac{2}{[n_r]([n_r]-1)}\sum_{1 \leq i < j \leq [n_r]}(f_{B(a,b)}(X_i, X_j) - \mu_{B(a,b)}),$$

with $g_{B(a,b)}(x) = \mathbb{E} f_{B(a,b)}(x, X_2)$, where $[n_r]$ stands for the integer part of $n_r = \mathbb{E} N_r = 2\lambda r$. Taking into account (3.3) for $d = 2, j = 1$ and (2.9), we get

$$\operatorname{Var} g_{B(a,b)}(X_1) = \sigma^{(1,2)}_{B(a,b)} - \mu^2_{B(a,b)}.$$

Since $N_r \xrightarrow[r \to \infty]{\mathbb{P}} \infty$ such that $N_r/[n_r] \xrightarrow[r \to \infty]{\mathbb{P}} 1$, Theorem VIII.7.1 in [26] tells us that in (5.7) the deterministic integer $[n_r]$ can be replaced by the random number $N_r$ without changing the limit $\mathcal{N}(0, \operatorname{Var} g_{B(a,b)}(X_1))$. Finally, a straightforward application of Slutsky's lemma completes the proof of Theorem 5.1. □

Notice that the asymptotic variance $\sigma^{(1,2)}_{B(a,b)} - \mu^2_{B(a,b)}$ of the number of intersection points under random centering in (5.6) is always smaller than the asymptotic variance $\sigma^{(1,2)}_{B(a,b)}$ obtained in (5.3) under deterministic centering.

Furthermore, instead of imposing conditions on the angles of the orientation vectors of the intersecting lines $\ell_{(p_i,\gamma_i)}$ and $\ell_{(p_j,\gamma_j)}$, in [24] the two angles at the intersecting points $x_i$ and $x_j$ of $\ell_{(p_i,\gamma_i)}$ and $\ell_{(p_j,\gamma_j)}$, respectively, with the $x$-axis are considered. More precisely, the ordered pair of angles $(\alpha, \beta)$ is considered, where $\alpha$ denotes the angle at $x_i \wedge x_j$ between the $x$-axis and the half-line from $\ell_{(p_i,\gamma_i)} \cap \ell_{(p_j,\gamma_j)}$ to the intersection point with the $x$-axis and, in analogy, $\beta$ is the angle at $x_i \vee x_j$.

**6. CLTs for Poisson–Voronoi tessellations.** In this section we give a brief overview of CLTs for Poisson–Voronoi tessellations (PVTs for short) in $\mathbb{R}^d$, that is, we consider a stationary PVT $\Phi$, generated by a stationary Poisson point process $\Psi$ in $\mathbb{R}^d$ with intensity $\lambda$. Notice that each cell of the Poisson–Voronoi tessellation $\Phi$ is defined as the closure of the set of all points in $\mathbb{R}^d$ which are closest to the nucleus of this cell, being a point of $\Psi$. Moreover, it can be shown that each cell is a convex ($d$-dimensional) polytope.

Let $\Psi_0$ denote the stationary (and isotropic) point process of vertices of the cells induced by $\Phi$ and let $\lambda_0$ denote its intensity, that is, $\lambda_0 = \mathbb{E}\Psi_0([0,1]^d)$. It is well known (see, e.g., [21]) that

$$(6.1) \qquad \lambda_0 = c_d \lambda \qquad \text{with } c_d = \frac{2^d \pi^{d-1}}{d+1} \frac{\kappa_{d^2}}{\kappa_{d^2-1}} \left(\frac{\kappa_{d-1}}{\kappa_d}\right)^d.$$

Consider a convex averaging sequence $W_n$ of sets in $\mathbb{R}^d$, that is, the sets $W_n$ are compact convex, increasing and contain a ball with unboundedly growing radius; see, for example, [6]. Since the $\beta$-mixing coefficient of a stationary PVT is exponentially decaying (see [9]), a CLT for the random number of vertices $\Psi_0(W_n)$ of a stationary PVT $\Phi$ within $W_n$ can be derived as $n \to \infty$. More precisely, in [12] it is shown that

$$(6.2) \quad \nu_d(W_n)^{-1/2}(\Psi_0(W_n) - \lambda_0 \nu_d(W_n)) \xrightarrow[n \to \infty]{\mathrm{d}} \mathcal{N}(0, \lambda_0(1 + c_d \sigma_d^2)),$$

where $\sigma_d^2$ is a constant only depending on the dimension $d$ which is expressible in terms of multiple integrals. A detailed discussion including numerical computations for the cases of $d=2$ and $d=3$ can be found in [10] and [12], respectively. The rounded values obtained there are $\sigma_2^2 = 0.5$ and $\sigma_3^2 = 5.084$, together with $c_2 = 2$ and $c_3 = 6.768$ from (6.1).

Notice that $\widehat{\lambda}_{0,n} = \Psi_0(W_n)/\nu_d(W_n)$ is an unbiased estimator for the intensity $\lambda_0$ and, by (6.2), $\widehat{\lambda}_{0,n}$ is asymptotically normally distributed with mean $\lambda_0$ and asymptotic variance $\lambda_0(1 + c_d\sigma_d^2)$. A simple transformation using (6.2) and $\widehat{\lambda}_{0,n} \xrightarrow[n\to\infty]{\mathbb{P}} \lambda_0$ yields

$$(6.3) \qquad 2\sqrt{\frac{\nu_d(W_n)}{1 + c_d\sigma_d^2}}(\sqrt{\widehat{\lambda}_{0,n}} - \sqrt{c_d\lambda}) \xrightarrow[n\to\infty]{\mathrm{d}} \mathcal{N}(0,1).$$

Similar to the confidence intervals derived in Corollary 4.2 for Poisson hyperplane processes, the CLT in (6.3) enables us to construct a $100(1-\alpha)\%$ (asymptotically exact) confidence interval for $\lambda$. Indeed, for any $\alpha \in (0,1)$, we have (for large enough $W_n$) that $\lambda \in I_n(\alpha)$ with probability $1-\alpha$, where $I_n(\alpha)c_d\nu_d(W_n)$ is given by

$$\left[\left(\sqrt{\Psi_0(W_n)} - \frac{z_{1-\alpha/2}}{2}\sqrt{1 + c_d\sigma_d^2}\right)^2, \left(\sqrt{\Psi_0(W_n)} + \frac{z_{1-\alpha/2}}{2}\sqrt{1 + c_d\sigma_d^2}\right)^2\right].$$

Notice that an alternative and general approach to CLTs in the Poisson–Voronoi case is given in [25]. Among other results, a CLT for a quite large class of functionals related to Poisson–Voronoi tessellations is obtained there.

**Acknowledgments.** We are grateful to the referees and an Associate Editor for their helpful comments and suggestions.

L. HEINRICH
INSTITUTE OF MATHEMATICS
UNIVERSITY OF AUGSBURG
D-86135 AUGSBURG
GERMANY
E-MAIL: lothar.heinrich@math.uni-augsburg.de

H. SCHMIDT
V. SCHMIDT
DEPARTMENT OF STOCHASTICS
UNIVERSITY OF ULM
D-89069 ULM
GERMANY
E-MAIL: hendrik.schmidt@uni-ulm.de
       volker.schmidt@uni-ulm.de